\documentclass[11pt,reqno]{amsart}

\usepackage{cite}
\usepackage{url}
\usepackage{amssymb}
\usepackage{amscd}
\usepackage{amsfonts}
\usepackage{amsmath}
\usepackage{eepic}
\usepackage{amsthm}
\usepackage{amsgen}
\usepackage[dvips]{color}
\usepackage{eucal}

 \DeclareMathOperator{\var}{\mathsf{var}} \DeclareMathOperator{\Id}{Eq}

\DeclareSymbolFont{rsfscript}{OMS}{rsfs}{m}{n}
\DeclareSymbolFontAlphabet{\mathrsfs}{rsfscript}

\numberwithin{equation}{section}

\def\bb{\mathbb}

\newtheorem*{Summary}{Theorem}
\newtheorem{Thm}{Theorem}[section]
\newtheorem{Prop}[Thm]{Proposition}
\newtheorem{Lemma}[Thm]{Lemma}

\newtheorem{Cor}[Thm]{Corollary}

\theoremstyle{remark}
\newtheorem{Rmk}{Remark}[section]
\newtheorem{Problem}{Problem}[section]

\def\om{\omega}

\def\cal{\mathcal}
\def\Ac{{\cal A}}
\def\Bc{{\cal B}}

\def\Hc{\mathrsfs{H}}

\def\Mc{{\cal M}}

\def\Rc{\mathrsfs{R}}
\def\Sc{{\cal S}}
\def\Tc{{\cal T}}
\def\Uc{{\cal U}}
\def\Vc{\mathbf{V}}

\def\si{\sigma}
\def\Si{\Sigma}

\def\la{\lambda}

\def\ol{\overline}

\def\H{\mathrm H}
\def\P{\mathrm P\!}
\def\he#1{#1#1^*}

\def\Rm{Rees matrix}
\def\sm{semi\-group}
\def\va{variet}

\def\fb{finitely based}
\def\ib{identity basis}
\def\nfb{non\-finitely based}

\def\TB{\ensuremath{\mathcal{T\kern-1pt B}_2^1}}
\def\TA{\ensuremath{\mathcal{T\kern-1pt A}_2^1}}
\def\BT{\ensuremath{\mathcal{B\kern-1pt T}}}
\def\BR{\ensuremath{\mathcal{B\kern-1pt R}}}
\def\BU{\ensuremath{\mathcal{B\kern-1pt U}}}
\def\FI{\ensuremath{\mathcal{FI}}}

\title[Matrix Identities Involving Multiplication and Transposition]%
{Matrix Identities Involving \\Multiplication and Transposition}

\author{K.~Auinger}
\address{Fakult\"at f\"ur Mathematik, Universit\"at Wien,
Nordbergstrasse 15,  A-1090 Wien, Austria}
\email{karl.auinger@univie.ac.at}
\author{I.~Dolinka}
\address{Department of Mathematics and Informatics, University of Novi Sad,
Trg Dositeja Obradovi\'ca 4, 21000 Novi Sad, Serbia}
\email{dockie@dmi.uns.ac.rs}
\author {M.~V.~Volkov}
\address{Faculty of Mathematics and Mechanics, Ural State University,
Lenina 51, 620083 Ekaterinburg, Russia}
\email{mikhail.volkov@usu.ru}

\begin{document}

\begin{abstract}
We study matrix identities involving multiplication and unary operations such as transposition
or Moore-Penrose inversion. We prove that in many cases such identities admit no finite basis.
\end{abstract}

\maketitle

\section*{Background and motivation}

Matrices\footnote{In this paper, the word `matrix' always means a square matrix of
finite size; also, in order to avoid trivialities, we always assume that the size
is at least two.} and matrix operations constitute basic tools for many branches of mathematics.
Important properties of matrix operations are often expressed in form of \emph{laws} or
\emph{identities} such as the associative law for multiplication of matrices. Studying
matrix identities that involve multiplication and addition is a classic research direction
that was motivated by several important problems in geometry and algebra (see \cite{Ami74}
for a survey of the origins of the theory) and that has eventually led to the profound and
beautiful theory of PI-rings \cite{DrFo04,GiZa05,KaRo05,Ro80}. Matrix identities involving
along with multiplication and addition also certain involution operations (such as taking
the usual or symplectic transpose of a matrix) have attracted much attention as well, see,
for instance, \cite{DaRa99,DaRa04,DrGi95,Gi90,Ro80}.

If one aims to classify matrix identities of a certain type, then a natural approach is
to look for a collection of `basic' identities such that all other identities would
follow from these basic identities. Such a collection is usually referred to as an
\emph{identity basis} or simply a \emph{basis}. For instance, all identities of matrices
over an infinite field involving only multiplication are known to follow from the associative
law, see \cite[Lemma~2]{GoMi78}. Thus, the associative law forms a basis of such `multiplicative'
identities. For identities involving both multiplication and addition, an explicit basis is known
for $2\times 2$-matrices (except the case of an infinite ground field of characteristic~2),
see \cite{Ra73,MaKu78,Dr81,Ko01,CoKo04}, and for $3\times 3$- and $4\times 4$-matrices over
a finite field, see \cite{Ge81,GeSi82}. However, for matrices of arbitrary size over a finite
field or a field of characteristic~0 the powerful results by Kruse-L'vov \cite{Kr73,Lv73}
and Kemer \cite{Ke87,Ke91} ensure at least the existence of a finite identity basis for such
identities.

In contrast, multiplicative identities of matrices over a finite field admit no finite basis.
This rather surprising fact was proved in the mid-1980s by the third author~\cite[Proposition~3]{V}
and Sapir~\cite[Corollary 6.2]{sapirburnside}. It is worth noting that methods used in~\cite{V}
and~\cite{sapirburnside} were very different but each of them sufficed to cover multiplicative
identities of matrices of every fixed size over every finite field.

In the present paper we study matrix identities involving multiplication and one or
two natural one-place operations such as taking various transposes or Moore-Penrose
inversion. For this we first have to adapt the methods of~\cite{V} and~\cite{sapirburnside}.
We present the corresponding results in Section~\ref{tools} while Section~\ref{preliminaries}
collects necessary preliminaries. Applications to the finite basis problem for matrix identities
are presented in Section~\ref{applications}. Both methods of Section~\ref{tools} are used here,
and it turns out that they in some sense complement one another since, in contrast to the case
of purely multiplicative identities, none of the methods alone suffices to cover, say, identities
of matrices of every size over every finite field involving multiplication and the usual
transposition of matrices.

Our main results may be a summarized as follows.

\begin{Summary}
Each of following sets of matrix identities admits no finite identity basis:
\begin{itemize}
\item the identities of $n\times n$-matrices over a finite field involving multiplication
and usual transposition;
\item the identities of $2n\times 2n$-matrices over a finite field involving multiplication
and symplectic transposition;
\item the identities of $2\times 2$-matrices over the field of complex numbers involving
either multiplication and Moore-Penrose inversion or multiplication, Moore-Penrose inversion
and Hermitian conjugation;
\item the identities of Boolean $n\times n$-matrices involving multiplication and transposition.
\end{itemize}
\end{Summary}

We mention in passing that tools developed in Section~\ref{tools} admit many further applications
that will be published in a separate paper.

As far as the theory of matrices is concerned, we use fairly standard concepts of linear algebra,
see \cite{Meyer}. Our proofs however involve some notions of equational logic and semigroup theory
which are briefly recalled in the next section.

\section{Preliminaries}\label{preliminaries}

\subsection{Vocabulary of equational logic}\label{vocabulary logic}

The concepts of an identity and an identity basis are intuitively clear. Nevertheless,
any precise reasoning about these concepts requires a formal framework, especially
when one aims at `negative' results as we do in this paper. Such a framework, provided
by equational logic, is concisely presented, e.g., in \cite[Chapter~II]{BuSa81}.
For the reader's convenience, we briefly overview the basic vocabulary of equational
logic in a form adapted to the use in this paper. Readers familiar with equational
logic may, of course, skip this overview.

A non-empty set $A$ endowed with operations
$f_1:\underbrace{A\times A\times\dots\times A}_{n_1}\to A$,\linebreak
$f_2:\underbrace{A\times A\times\dots\times A}_{n_2}\to A$, \dots\
is called an \emph{algebraic structure of type $(n_1,n_2,\dots)$ with the carrier $A$}.
Algebraic structures considered in this paper are mostly of types $(2,1)$ or $(2,1,1)$ so
that they have one binary operation and one or two unary operations. Since binary
operations involved are always assumed to satisfy the associative law, our
structures are \emph{semigroups} equipped with one or two additional unary
operations; for brevity, we shall call such structures \emph{unary semigroups}.

We formally define notions related to unary semigroup identities
only for the case of one unary operation leaving the straightforward
modification for the case of two unary operations to the reader.
Given a countably infinite set $X$, we define the set $T(X)$ of all
\emph{unary semigroup terms} over $X$ as follows:
\begin{itemize}
\item every $x\in X$ is a unary semigroup term and so is $x^*$;
\item if $u$ and $v$ are unary semigroup terms,
then so is their concatenation $uv$;
\item if $u$ is a unary semigroup term, then so is
$(u)^*$.
\end{itemize}
The algebraic structure $\mathcal{T}(X)=\langle T(X),\cdot,{}^*\rangle$
of type $(2,1)$ whose binary operation $\cdot$ is concatenation and whose
unary operation is $u\mapsto(u)^*$ is called the \emph{free unary \sm} over $X$.
By a \emph{unary \sm\ identity} over $X$ we mean a formal expression $u=v$
where $u,v\in T(X)$. A unary \sm\ $\mathcal{S}=\langle S,\cdot,{}^*\rangle$
\emph{satisfies} the identity $u=v$ if the equality $\varphi(u)=\varphi(v)$
holds in $\mathcal{S}$ under all possible homomorphisms $\varphi:\mathcal{T}(X)
\to\mathcal{S}$. Given $\mathcal{S}$, we denote by $\Id\mathcal{S}$ the set
of all unary \sm\ identities it satisfies.

Given any collection $\Sigma$ of unary \sm\ identities, we say that an identity $u=v$
\emph{follows} from $\Sigma$ or that $\Sigma$ \emph{implies} $u=v$ if every
unary \sm\ satisfying all identities of $\Sigma$ satisfies the identity $u=v$ as
well. Birkhoff's completeness theorem of equational logic (see \cite[Theorem 14.17]{BuSa81})
shows that this notion (which we have given a semantic definition) can be captured by
a very transparent set of inference rules. These rules in fact formalize the most
natural things one does with identities: substitution of a term for a variable,
application of operations to identities (such as, say, multiplying an identity
through on the right by the same term) and using symmetry and transitivity of equality.
We need not going into more detail here because the completeness theorem
is not utilized in this paper.

Given a unary \sm\ $\mathcal{S}$, an \emph{\ib} for $\mathcal{S}$ is any set
$\Sigma\subseteq\Id\mathcal{S}$ such that every identity of $\Id\mathcal{S}$
follows from $\Sigma$. A unary \sm\ $\mathcal{S}$ is said to be \emph{\fb} if it possesses
a finite \ib; otherwise $\Sc$ is called \emph{\nfb}.

The class of all unary semigroups satisfying all identities from a given set
$\Sigma$ of unary \sm\ identities is called the \emph{variety defined by $\Sigma$}.
A variety is is said to be \emph{\fb} if it can be defined by a finite set
of identities; otherwise it is called \emph{\nfb}.

It is easy to see that the satisfaction of an identity is inherited by forming
direct products, taking unary subsemigroups and homomorphic images so that
each variety is closed under these operators. In fact, varieties can be characterized
by this closure property (the HSP-theorem, see \cite[Theorem 11.9]{BuSa81}).
Given a unary \sm\ $\mathcal{S}$, the variety defined by $\Id\Sc$
is the \emph{variety generated by $\Sc$}; we denote this variety by
$\var\Sc$. From the HSP-theorem it follows that every member of $\var\Sc$
is a homomorphic image of a unary subsemigroup of a direct
product of several copies of $\Sc$. Observe also that a unary \sm\ and
the variety it generates are simultaneously finitely or \nfb.

A variety $\Vc$ is said to be \emph{locally finite} if every finitely
generated member of $\Vc$ is finite. A  finite unary \sm\ is called
\emph{inherently \nfb} if it is not contained in any finitely based
locally finite variety. Since the variety generated by a finite unary
\sm\ is locally finite (this is an easy consequence of the HSP-theorem,
see \cite[Theorem 10.16]{BuSa81}), the property of being inherently \nfb\
implies the property of being \nfb; in fact, the former property is much stronger.

\subsection{Vocabulary of semigroup theory}\label{vocabulary semigroup}

Since the algebraic structures we deal with are semigroups with some
additional operation(s), we freely use the standard terminology and
notation of semigroup theory, mostly following the early chapters of
the textbook~\cite{CP}. It should be noted, however, that our presentation
is to a reasonable extent self-contained and does not require any specific
\sm-theoretic background.

In general, considering a unary semigroup $\mathcal{S}=\langle S,\cdot,{}^*\rangle$,
we do not assume any additional identities involving the unary operation ${}^*$.
If the identities $(xy)^* = y^*x^*$ and $(x^*)^* = x$ happen to hold in $\mathcal{S}$,
in other words, if the unary operation $x\mapsto x^*$ is an involutory anti-automorphism
of the semigroup $\langle S,\cdot\rangle$, we call $\mathcal{S}$ an \emph{involutory semigroup}.
If, in addition, the identity $x=xx^*x$ holds, $\mathcal{S}$ is said to be a \emph{regular
$*$-semigroup}. Each group, subject to its inverse operation $x\mapsto x^{-1}$, is
an involutory semigroup, even a regular $*$-semigroup; throughout the paper,
any group is considered as a unary semigroup with respect to this inverse unary
operation.

A wealth of examples of involutory semigroups and regular
$*$-semigroups can be obtained via the following `unary' version
of the well known \Rm\ construction\footnote{See \cite[Section~3.1]{CP}
for a description of the construction in the plain semigroup
case and \cite[Section~2]{GePe85} for a detailed analysis of its unary
version.}. Let $\mathcal{G}=\langle G,\cdot,{}^{-1}\rangle$ be a
group, $0$ a symbol beyond $G$, and $I$ a non-empty set. We
formally set $0^{-1}=0$. Given an $I\times I$-matrix $P=(p_{ij})$
over $G\cup\{0\}$ such that $p_{ij}=p_{ji}^{-1}$ for all $i,j\in
I$, we define a multiplication~$\cdot$ and a unary operation
${}^*$ on the set $(I\times G\times I)\cup\{0\}$ by the following
rules:
$$a\cdot 0=0\cdot a=0\ \text{ for all $a\in (I\times G\times I)\cup \{0\}$},$$
\begin{gather*}
(i,g,j)\cdot(k,h,\ell)=\left\{\begin{array}{cl}
(i,gp_{jk}h,\ell)&\ \text{if}\ p_{jk}\ne0,\\
0 &\ \text{if}\ p_{jk}=0;
\end{array}\right.\\
(i,g,j)^* = (j,g^{-1},i),\ 0^* = 0.
\end{gather*}
It can be easily checked that $\langle(I\times G\times I)\cup
\{0\},\cdot,{}^*\rangle$ becomes  an involutory semigroup; it will
be a regular $*$-semigroup precisely when $p_{ii}=e$ (the identity
element of the group $\mathcal{G}$) for all $i\in I$. We denote
this unary semigroup  by $\Mc^0(I,\mathcal{G},I;P)$ and call it
the \emph{unary \Rm\ \sm\ over $\mathcal{G}$ with the sandwich
matrix $P$}. If the involved group $\mathcal G$ happens to be the
trivial group $\mathcal{E}=\{e\}$, then we shall ignore the
group entry and represent the non-zero elements of such a Rees
matrix semigroup by the pairs $(i,j)$ with $i,j\in I$.

In this paper, the 10-element unary \Rm\ semigroup over the
trivial group $\mathcal{E}=\{e\}$ with the sandwich matrix
$$\begin{pmatrix}
e & e & e\\
e & e & 0\\
e & 0 & e
\end{pmatrix}$$
plays a key role; we denote this \sm\ by $\mathcal{K}_3$. Thus,
subject to the convention mentioned above, $\mathcal{K}_3$
consists of the nine pairs $(i,j)$, $i,j\in\{1,2,3\}$, and the
element $0$, and the operations restricted to its non-zero
elements can be described as follows:
\begin{gather}
\label{operations in C3}
(i,j)\cdot(k,\ell)=\left\{\begin{array}{cl}
(i,\ell)&\ \text{if}\ (j,k)\ne(2,3),(3,2),\\
0 &\ \text{otherwise};
\end{array}\right.\\
(i,j)^* = (j,i)\notag.
\end{gather}

Another unary semigroup that will be quite useful in the sequel is
the \emph{free involutory semigroup} $\FI(X)$ on a given
alphabet $X$. It can be constructed as follows.  Let
$\overline{X}=\{x^*\mid x\in X\}$ be a disjoint copy of $X$ and
define $(x^*)^*=x$ for all $x^*\in \overline{X}$. Then
$\FI(X)$ is the free semigroup $(X\cup\overline{X})^+$
endowed with an involution ${}^*$ defined by
$$(x_1\cdots x_m)^* = x_m^*\cdots x_1^*$$
for all $x_1,\dots,x_m\in X\cup \overline{X}$. See
\cite[Section~3]{GePe85} for more details on $\FI(X)$.

We will refer to elements of $\FI(X)$ as \emph{involutory words over
$X$} while elements of the free semigroup $X^+$ will be referred
to as (plain semigroup) words over $X$.

\subsection{A property of matrices of rank 1}
\label{rank 1 lemma}

Given a field $\mathcal{K}=\langle K,+,\cdot\rangle$, we
denote the set of all $n\times n$-matrices over $\mathcal{K}$ by
$\mathrm{M}_n(\mathcal{K})$. As mentioned in the introduction,
in order to avoid trivialities, we always assume that $n\ge 2$.

We conclude our preparations with registering a simple property of rank~1
matrices. This property is, of course, known, but we do provide a proof
for the sake of completeness.
\begin{Lemma}
\label{rank 1} If a matrix $A\in\mathrm{M}_n(\mathcal{K})$ has
rank $1$, then $A^2BA=ABA^2$ for any matrix $B\in\mathrm{M}_n(\mathcal{K})$.
\end{Lemma}

\begin{proof}
Consider the natural action of $A$ on the vector space $K^n$. The null-space $N(A)=\{x\in K^n\mid xA=0\}$ has dimension
$n-1$ whence the generalized eigenspace of $A$ corresponding to the eigenvalue $0$ coincides with either $K^n$ or
$N(A)$. In the former case $A^2=0$ and $A^2BA=0=ABA^2$ for any matrix $B$. In the latter case $K^n$ decomposes into the
direct sum of $N(A)$ and the range $R(A)=\{xA\mid x\in K^n\}$, see \cite[Section~5.10]{Meyer}. Then $R(A)$ is in fact
the (generalized) eigenspace of $A$ corresponding to a non-zero eigenvalue $\alpha\in K$ and the minimal polynomial of
$A$ is $x(x-\alpha)$. Thus, $A$ satisfies the equation $A^2-\alpha A=0$, whence $A^2BA=\alpha ABA=AB\alpha A=ABA^2$, as
required.
\end{proof}

Let $\mathrm{L}_n(\mathcal{K})$ denote the set of all $n\times
n$-matrices of rank at most~1 over $\mathcal{K}$. Adding the
identity matrix to $\mathrm{L}_n(\mathcal{K})$ we get a set
which we denote by $\mathrm{L}^1_n(\mathcal{K})$. Clearly, it is
closed under matrix multiplication. From Lemma~\ref{rank 1} we
immediately obtain
\begin{Cor}
\label{identity for rank 1} The semigroup
$\langle\mathrm{L}^1_n(\mathcal{K}),\cdot\rangle$ satisfies the
identity
\begin{equation}
\label{abelian groups} x^2yx=xyx^2.
\end{equation}
\end{Cor}
Observe that every group satisfying \eqref{abelian groups} is
abelian.

\section{Tools}\label{tools}

\subsection{A unary version of the critical semigroup method}
Here we present a `unary' modification of the approach used in~\cite{V}. According
to the classification proposed in the survey~\cite{volkovjaponicae}, this approach
is referred to as the \emph{critical semigroup method}.

The formulation of the corresponding result involves two simple
operators on unary semigroup varieties. For any unary semigroup
$\mathcal{S}=\langle S,\cdot,{}^*\rangle$ we denote by
$\H(\mathcal{S})$ the unary subsemigroup of $\mathcal{S}$ which is
generated by all elements of the form $xx^*$, where $x\in S$. We
call $\H(\mathcal{S})$ the \emph{Hermitian subsemigroup} of
$\mathcal{S}$. For any variety $\Vc$ of unary semigroups, let
$\H(\Vc)$ be the subvariety of $\Vc$ generated by all Hermitian
subsemigroups of members of $\Vc$. Likewise, given a positive
integer $n$, let $\P_n(\mathcal{S})$ be the unary subsemigroup of
$\mathcal{S}$ which is generated by all elements of the form
$x^n$, where $x\in S$, and let $\P_n(\Vc)$ be the subvariety of
$\Vc$ generated by all subsemigroups $\P_n(\mathcal{S})$, where
$\mathcal{S}\in \Vc$.

The following easy observation will be useful in
the sequel as it helps calculating the effect of the operators
$\H$ and $\P_n$.

\begin{Lemma}\label{Lemma 3.1}
$\H(\var\Sc)=\var\H(\Sc)$ and $\P_n(\var\Sc)= \var\P_n(\Sc)$ for
every unary semigroup $\Sc$ and for each $n$.
\end{Lemma}

\begin{proof}
The non-trivial part of the first claim is the inclusion
$\H(\var\Sc)\subseteq\var\H(\Sc)$. Let $\Tc\in\var\Sc$, then $\Tc$
is a homomorphic image of a unary subsemigroup $\Uc$ of a direct
product of several copies of $\Sc$. But then $\H(\Tc)$ is a
homomorphic image of $\H(\Uc)$. As is easy to see, $\H(\Uc)$ is a
unary subsemigroup of a direct product of several copies of
$\H(\Sc)$. Thus $\H(\Tc)\in\var\H(\Sc)$. Since this holds for an
arbitrary $\Tc\in\var\Sc$, we conclude that
$\H(\var\Sc)\subseteq\var\H(\Sc)$. The second assertion can be
treated in a completely similar way.
\end{proof}

We are now ready to state the main result of this subsection.

\begin{Thm}\label{Theorem 2.1}
Let $\Vc$ be any unary semigroup variety such that
$\mathcal{K}_3\in\Vc$. If either
\begin{itemize}
\item there exists a group $\mathcal{G}$  such that $\mathcal{G}\in\Vc$ but
$\mathcal{G}\notin\H(\Vc)$

\hbox to 10 pt {or}

\item there exist a positive integer $d$ and a group $\mathcal{G}$ of exponent
dividing $d$ such that $\mathcal{G}\in\Vc$ but
$\mathcal{G}\notin\P_d(\Vc)$,
\end{itemize}
then $\Vc$ has no finite basis of identities.
\end{Thm}

\begin{proof}
Assume first that there exists a group $\mathcal{G}\in\Vc$ for
which $\mathcal{G}\notin\H(\Vc)$.

\medskip

1. First we recall the basic idea of `the critical semigroup
method' in the unary setting. Suppose that ${\Vc}$ is \fb. If
$\Si$ is a finite identity basis of the variety ${\Vc}$ then there
exists a positive integer $\ell$ such that all identities from
$\Si$ depend on at most $\ell$ letters. Therefore identities from
$\Si$ hold in a unary semigroup $\mathcal{S}$ whenever all
$\ell$-generated unary subsemigroups of $\mathcal{S}$ satisfy
$\Si$. In other words, $\mathcal{S}$ belongs to $\Vc$ whenever all
of its $\ell$-generated unary subsemigroups are in ${\Vc}$. We see
that in order to prove our theorem it is sufficient to construct,
for any given positive integer $k$, a unary \sm\
$\mathcal{T}_k\notin{\Vc}$ for which all $k$-generated unary
sub\sm s of $\mathcal{T}_k$ belong to ${\Vc}$.

\medskip

 2. Fix an identity $u(x_1,\ldots,x_m) = v(x_1,\ldots,x_m)$ that holds in ${\H(\Vc)}$
but fails in the group $\mathcal{G}=\langle G,\cdot,{}^{-1}\rangle$. The latter means
that, for some $g_1,\dots,g_m\in G$, substitution of $g_i$ for $x_i$ yields
\begin{equation}\label{2.1}
u(g_1,\ldots,g_m) \ne v(g_1,\ldots,g_m).
\end{equation}
Now, for each positive integer $k$, let $n =\max\{4,2k+1\}$,
$I=\{1,\dots,nm\}$ and consider the unary \Rm\ \sm \
$\mathcal{T}_k=\Mc^0(I,\mathcal{G},I;P_k)$ over the group
$\mathcal{G}$ with the sandwich matrix
$$P_k=\left(\begin{array}{ccccccc}
M_n(g_1) & E_n & O_n & O_n & \cdots & O_n & E_n^T \\
E_n^T & M_n(g_2) & E_n & O_n & \cdots & O_n & O_n \\
O_n & E_n^T & M_n(g_3) & E_n & \cdots & O_n & O_n \\
\vdots & \vdots & \vdots & \vdots & \ddots & \vdots & \vdots \\
O_n & O_n & O_n & O_n & \cdots & E_n & O_n\\
O_n & O_n & O_n & O_n & \cdots & M_n(g_{m-1}) & E_n \\
E_n & O_n & O_n & O_n & \cdots & E_n^T & M_n(g_m)
\end{array}\right),$$
where $O_n$ is the zero $n\times n$-matrix, $E_n$ is the $n\times
n$-matrix having $e$ (the identity of $\mathcal{G}$) in the
position $(n,1)$ and 0 in all other positions, $E_n^T$ is the
transpose of $E_n$, and $M_n(g)$ denotes the $n\times n$-matrix of
the form
$$M_n(g)=\left(\begin{array}{ccccccc}
e & g &  0 & \cdots & 0 & 0 & e \\
g^{-1} & e & e & \cdots & 0 & 0 & 0 \\
0 & e & e & \cdots & 0 & 0 & 0 \\
\vdots & \vdots & \vdots & \ddots & \vdots & \vdots & \vdots \\
0 & 0 & 0 & \cdots & e & e & 0 \\
0 & 0 & 0 & \cdots & e & e & e \\
e & 0 & 0 & \cdots & 0 & e & e
\end{array} \right).$$
(This construction is in a sense a combination of those of the
first and the third authors' papers \cite{A1} and \cite{V}.) We
are going to prove that $\mathcal{T}_k$ enjoys the two properties
needed, namely, it does not belong to ${\Vc}$, but each
$k$-generated unary sub\sm \ of $\mathcal{T}_k$ lies in ${\Vc}$.

\medskip

3. In order to prove that $\mathcal{T}_k\notin\Vc$, we construct
an identity that holds in ${\Vc}$, but fails in $\mathcal{T}_k$.
Consider the following $m$ terms in $mn$ letters
$x_1,\dots,x_{mn}$ (the square brackets in these terms
serve only to improve readability):

\smallskip

\leftline{$w_1=[\he{x_1}\cdots\he{x_n}] [\he{(x_{n+1}\cdots
x_{2n})}\he{x_{2n}}]\cdots$} \rightline{[$\he{(x_{(m-1)n+1}\cdots
x_{mn})}\he{x_{mn}}]$,}

\leftline{$w_2=[\he{(x_1\cdots x_n)}\he{x_n}]
[\he{x_{n+1}}\cdots\he{x_{2n}}]\times$}
\centerline{$[\he{(x_{2n+1}\cdots x_{3n})}\he{x_{3n}}]\cdots$}
\rightline{$[\he{(x_{(m-1)n+1}\cdots x_{mn})}\he{x_{mn}}]$,}

\leftline{$w_3=[\he{(x_1\cdots x_n)}\he{x_n}] [\he{(x_{n+1}\cdots
x_{2n})}\he{x_{2n}}]\times$}
\rightline{$[\he{x_{2n+1}}\cdots\he{x_{3n}}]\cdots[\he{(x_{(m-1)n+1}\ldots
x_{mn})}\he{x_{mn}}]$,}

\hbox to 5.0 in{\dotfill }

\leftline{$w_m=[\he{(x_1\cdots x_n)}\he{x_n}] [\he{(x_{n+1}\cdots
x_{2n})}\he{x_{2n}}]\ldots$} \centerline{$[\he{(x_{(m-2)n+1}\cdots
x_{(m-1)n})} \he{x_{(m-1)n}}]\times$}
\rightline{$[\he{x_{(m-1)n+1}}\cdots\he{x_{mn}}]$.}

\smallskip

\noindent Substituting $w_i$ for $x_i$ in $u$ respectively $v$, we
get the identity
\begin{equation}\label{2.2}
u(w_1,\ldots,w_m) = v(w_1,\ldots,w_m) \end{equation} which holds
in the \va y ${\Vc}$. Indeed, if we take any
$\mathcal{S}\in{\Vc}$, then, since $\he{s}\in\H(\mathcal{S})$ for
any $s$ in $\mathcal{S}$, all the values of $w_i$ belong to the
Hermitian sub\sm\ $\H(\mathcal{S})$ of $\mathcal{S}$. This sub\sm,
however, lies in ${\H(\Vc)}$, and therefore, satisfies the identity
$u = v$.

Now we shall show that \eqref{2.2} fails in $\mathcal{T}_k$.
Indeed, substituting $(i,e,i)\in\mathcal{T}_k$ for $x_i$, we
calculate that the value of every term of the form
$$\he{(x_{(j-1)n+1}\cdots x_{jn})}\he{x_{jn}}$$
is equal to $((j-1)n+1,e,jn)$ while the value of each term of the
form
$$\he{x_{(j-1)n+1}}\cdots \he{x_{jn}}$$
is equal to  $((j-1)n+1,g_j,jn)$. Hence the value of $w_j$ is just
$(1,g_j,mn)$. Therefore, under this substitution, the left hand
part of (\ref{2.2}) takes the value $(s,u(g_1,\ldots,g_m),t)$ for
suitable $s,t\in \{1,mn\}$ while the value of the right hand part
of (\ref{2.2}) is $(s',v(g_1,\ldots,g_m),t')$ (again for suitable
$s',t'\in\{1,mn\}$). In view of the inequality (\ref{2.1}), these
elements do not coincide in $\mathcal{T}_k$.

\smallskip

4. It remains to prove that each $k$-generated unary sub\sm\ of
$\mathcal{T}_k$ lies in ${\Vc}$. For every $m$-tuple
$(\la_1,\ldots,\la_m)$ of positive integers satisfying
\begin{equation} \label{2.3}
1\leq\la_1\leq n<\la_2\leq 2n<\la_3\leq\ldots (m-1)n <\la_m\leq
mn,
\end{equation}
consider the unary sub\sm \ $\mathcal{T}_k(\la_1,\ldots,\la_m)$ of
$\mathcal{T}_k$ consisting of 0 and all triples $(i,g,j)$ such
that $g\in G$ and $i,j\notin\{\la_1,\ldots,\la_m\}$.
Using that $2k<n$ according to our choice of $n$, one concludes
that any given $k$ elements of $\mathcal{T}_k$ must be contained
in $\mathcal{T}_k(\la_1,\ldots,\la_m)$ for suitable
$\la_1,\ldots,\la_m$. Thus it is sufficient to prove that each
semigroup  of the form $\mathcal{T}_k(\la_1,\ldots,\la_m)$ belongs
to the variety ${\Vc}$.

Let us fix positive integers $\la_1,\ldots,\la_m$ satisfying
(\ref{2.3}). When multiplying triples from
$\mathcal{T}_k(\la_1,\ldots,\la_m)$, the
$\la_1^{\mathrm{th}},\dots,\la_m^{\mathrm{th}}$ rows and columns
of the sandwich matrix $P_k$ are never involved. Therefore we can
identify $\mathcal{T}_k(\la_1,\ldots,\la_m)$ with the unary \Rm\
\sm\ $\Mc^0(I',\mathcal{G},I';P'_k)$ over the group $\mathcal{G}$
where $I'=I\setminus\{\la_1,\dots, \la_m\}$ and the sandwich
matrix $P'_k= P_k(\la_1,\ldots,\la_m)$ is obtained from $P_k$ by
deleting its $\la_1^{\mathrm{th}},\dots,\la_m^{\mathrm{th}}$ rows
and columns. Note that by \eqref{2.3} exactly one row and one
column of each block $M_n(g_i)$ is deleted.

Now we transform the matrix $P_k(\la_1,\ldots,\la_m)$ as follows.
For each $i$ such that $(i-1)n+2<\la_i$, we multiply successively
\begin{align}
\label{transform}
&\text{the row $((i-1)n+2)$ by $g_i$ from the left and}\notag\\[-1ex]
&\text{the column $((i-1)n+2)$ by $g_i^{-1}$ from the right;}\notag\\[-.5ex]
&\text{the row $((i-1)n+3)$ by $g_i$ from the left and}\notag\\[-1ex]
&\text{the column $((i-1)n+3)$ by $g_i^{-1}$ from the right;}\\[-.5ex]
&\hbox to 3.0 in{\dotfill }\notag\\[-.5ex]
&\text{the row $(\la_i-1)$ by $g_i$ from the left and}\notag\\[-1ex]
&\text{the column $(\la_i-1)$ by $g_i^{-1}$ from the right.}\notag
\end{align}

In order to help the reader to understand the effect of the
transformations~\eqref{transform}, we illustrate their action on
the block obtained from $M_n(g_i)$ by removing the
$\la_i^{\mathrm{th}}$ row and column in the following scheme in
which $\lambda_i$ has been chosen to be equal to $(i-1)n+5$. (The
transformations have no effect beyond $M_n(g_i)$ because all the
rows and columns of $P_k(\la_1,\ldots,\la_m)$ involved
in~\eqref{transform} have non-zero entries only within
$M_n(g_i)$.)

{\small
\begin{center}
\begin{tabular}{cc}
The block obtained from $M_n(g_i)$ by erasing  & After the first\\
the $((i-1)n+5)^{\mathrm{th}}$ row and column  & transformation \\[1ex]
$\begin{pmatrix}
e        & g_i & 0 & 0 &  0 &\cdots & 0 & e \\
g_i^{-1} & e   & e & 0 &  0 &\cdots & 0 & 0 \\
0        & e   & e & e &  0 &\cdots & 0 & 0 \\
0        & 0   & e & e &  0 &\cdots & 0 & 0 \\
0        & 0   & 0 & 0 &  e &\cdots & 0 & 0 \\
\vdots   & \vdots & \vdots & \vdots & \vdots &\ddots & \vdots & \vdots \\
0        & 0   & 0 & 0 &  0 &\cdots & e & e \\
e        & 0   & 0 & 0 &  0 &\cdots & e & e
\end{pmatrix}$
& $\begin{pmatrix}
e        & e   & 0 & 0 &  0 &\cdots & 0 & e \\
e   & e  & g_i & 0 &  0 &\cdots & 0 & 0 \\
0   & g_i^{-1} & e & e &  0 &\cdots & 0 & 0 \\
0        & 0   & e & e &  0 &\cdots & 0 & 0 \\
0        & 0   & 0 & 0 &  e &\cdots & 0 & 0 \\
\vdots   & \vdots & \vdots & \vdots & \vdots &\ddots & \vdots & \vdots \\
0        & 0   & 0 & 0 &  0 &\cdots & e & e \\
e        & 0   & 0 & 0 &  0 &\cdots & e & e
\end{pmatrix}$\\
&\\
After the second & After the third\\
transformation   & transformation\\[1ex]
$\begin{pmatrix}
e        & e   & 0 & 0 &  0 &\cdots & 0 & e \\
e        & e   & e & 0 &  0 &\cdots & 0 & 0 \\
0        & e & e & g_i &  0 &\cdots & 0 & 0 \\
0 & 0   & g_i^{-1} & e &  0 &\cdots & 0 & 0 \\
0        & 0   & 0 & 0 &  e &\cdots & 0 & 0 \\
\vdots   & \vdots & \vdots & \vdots & \vdots &\ddots & \vdots & \vdots \\
0        & 0   & 0 & 0 &  0 &\cdots & e & e \\
e        & 0   & 0 & 0 &  0 &\cdots & e & e
\end{pmatrix}$
& $\begin{pmatrix}
e        & e   & 0 & 0 &  0 &\cdots & 0 & e \\
e        & e   & e & 0 &  0 &\cdots & 0 & 0 \\
0        & e   & e & e &  0 &\cdots & 0 & 0 \\
0        & 0   & e & e &  0 &\cdots & 0 & 0 \\
0        & 0   & 0 & 0 &  e &\cdots & 0 & 0 \\
\vdots   & \vdots & \vdots & \vdots & \vdots &\ddots & \vdots & \vdots \\
0        & 0   & 0 & 0 &  0 &\cdots & e & e \\
e        & 0   & 0 & 0 &  0 &\cdots & e & e
\end{pmatrix}$
\end{tabular}
\end{center}

} \noindent Now it should be clear that also in general the
transformations~\eqref{transform} result in a matrix $Q_k$ all of
whose non-zero entries are equal to $e$. On the other hand, it is
known (see, e.\,g., \cite[Proposition~6.2]{A1}) that the
transformations~\eqref{transform} of the sandwich matrix do not
change the unary \sm\ $\mathcal{T}_k(\la_1,\ldots,\la_m)$; in
other words $\mathcal{T}_k(\la_1,\ldots,\la_m)$ is isomorphic to
the $I'\times I'$ unary \Rm \ \sm \ $\mathcal{R}_k$ over
$\mathcal{G}$ with the sandwich matrix $Q_k$. Let $\mathcal{U}_k$
be the $I'\times I'$ unary \Rm \ \sm\ over the trivial group
$\mathcal{E}$ with the sandwich matrix $Q_k$. It is easy to check
that the mapping $\mathcal{G}\times \mathcal{U}_k\rightarrow
\mathcal{R}_k$ defined by
$$(g,(i,j))\mapsto (i,g,j),\quad  (g,0)\mapsto 0$$
for all $g\in G$, $i,j\in I'$, is a unary semigroup
homomorphism onto $\mathcal{R}_k$. Now we note that $\mathcal{G}
\in {\Vc}$ and $\mathcal{U}_k$ belongs to the \va y generated by
$\mathcal{K}_3$ (see \cite[Theorem 5.2]{A1}). This yields
$\mathcal{T}_k(\la_1,\ldots,\la_m)\cong\mathcal{R}_k\in{\Vc}$.

\smallskip

The case when there exist a positive integer $d$ and a group
$\mathcal{G}$ of exponent dividing $d$ such that
$\mathcal{G}\in\Vc$ but $\mathcal{G}\notin\P_d(\Vc)$ can be
treated in a very similar way. The construction of the critical
semigroups remains the same, and the only modification to be made
in the rest of the proof is to replace the terms $w_i$ above by
the following plain semigroup words:

\smallskip

\leftline{$w_1=[x_1^d\cdots x_n^d][x_{n+1}\cdots x_{2n}]^d\cdots
[x_{(m-1)n+1}\cdots x_{mn}]^d$,} \leftline{$w_2=[x_1\cdots
x_n]^d[x_{n+1}^d\cdots x_{2n}^d]\cdots [x_{(m-1)n+1}\cdots
x_{mn}]^d$,}

\hbox to 3.3 in {\dotfill }

\leftline{$w_m=[x_1\cdots x_n]^d[x_{n+1}\cdots
x_{2n}]^d\cdots[x_{(m-1)n+1}^d\cdots x_{mn}^d]$.}

\smallskip

\noindent (These words already have been used in the plain
semigroup case by the third author~\cite{V}.)
\end{proof}

\subsection{A unary version of the method of inherently
nonfinitely based semigroups} Here we prove a sufficient condition for an involutory
semigroup to be inherently \nfb\ and exhibit two concrete examples of involutory semigroups
satisfying this condition. These examples will be essential in our applications in Section~\ref{applications}.

Let $x_1,x_2,\dots,x_n,\dots$ be a sequence of letters. The
sequence $\{Z_n\}_{n=1,2,\dots}$ of \emph{Zimin words} is defined
inductively by $Z_1=x_1$, $Z_{n+1}=Z_nx_{n+1}Z_n$. We say that an
involutory word $v$ is an \emph{involutory isoterm for a unary
semigroup $\mathcal{S}$} if the only involutory word $v'$ such
that $\mathcal{S}$ satisfies the involutory semigroup identity
$v=v'$ is the word $v$ itself.

\begin{Thm}
\label{Theorem 2.2} Let $\mathcal{S}$ be a finite involutory
semigroup. If all Zimin words are involutory isoterms for
$\mathcal{S}$, then $\mathcal{S}$ is inherently \nfb.
\end{Thm}

\begin{proof}
Arguing by contradiction, suppose that $\mathcal{S}$ belongs to a
\fb\ locally finite variety $\Vc$. If $\Si$ is a finite identity
basis of $\Vc$, then there exists a positive integer $\ell$ such
that all identities from $\Si$ depend on at most $\ell$ letters.
Clearly, all identities in $\Si$  hold in $\mathcal{S}$. In the
following, our aim will be to construct, for any given positive
integer $k$, an infinite, finitely generated involutory \sm\
$\mathcal{T}_k$ which satisfies all identities in at most $k$
variables that hold in $\mathcal{S}$. In particular,
$\mathcal{T}_\ell$ will satisfy all identities from $\Si$. This
yields a contradiction, as then we must conclude that
$\mathcal{T}_\ell\in\Vc$, which is impossible by the local
finiteness of $\Vc$.

We shall employ a construction invented by
Sapir~\cite{sapirburnside}, see also his lecture
notes~\cite{Sapir}. We fix $k$ and let $r=6k+2$. Consider the
$r^2\times r$-matrix $M$ shown in the left hand part
of Fig.~\ref{matrices}.
\begin{figure}[b]
$$M=\begin{pmatrix}
  1 & 1 & \cdots & 1 & 1\\

  \vdots & \vdots & \ddots & \vdots & \vdots \\
  1 & r & \cdots & 1 & r\\
  2 & 1 & \cdots & 2 & 1\\

  \vdots & \vdots & \ddots & \vdots & \vdots\\
  2 & r & \cdots & 2 & r\\
  \vdots & \vdots & \ddots & \vdots & \vdots\\
  r & 1 & \cdots & r & 1\\

  \vdots & \vdots & \ddots & \vdots & \vdots\\
  r & r & \cdots & r & r\\
\end{pmatrix} \qquad
 M_{A}=\begin{pmatrix}
  a_{11} & a_{12} & \cdots & a_{1r-1} & a_{1 r} \\

  \vdots & \vdots & \ddots & \vdots & \vdots \\
  a_{11} & a_{r2} & \cdots & a_{1r-1} & a_{rr}\\
  a_{21} & a_{12} & \cdots & a_{2r-1} & a_{1r}\\

  \vdots & \vdots & \ddots & \vdots & \vdots \\
  a_{21} & a_{r2} & \cdots & a_{2r-1} & a_{rr}\\
  \vdots & \vdots & \ddots & \vdots & \vdots \\
  a_{r1} & a_{12} & \cdots & a_{rr-1} & a_{1r}\\

  \vdots & \vdots & \ddots & \vdots & \vdots \\
  a_{r1} & a_{r2} & \cdots & a_{rr-1} & a_{rr}\\
\end{pmatrix}$$
\caption{The matrices $M$ and $M_A$}\label{matrices}
\end{figure}
All odd columns of $M$ are identical and equal to the transpose of
the row $(1,1,\ldots,1,2,2,\ldots,2,\ldots,r,r,\ldots,r)$ where
each number occurs $r$ times. All even columns of $M$ are
identical and equal to the transpose of the row
$(1,2,\ldots,r,1,2,\ldots,r,\ldots,1,2,\ldots,r)$ in which the
block $1,2,\dots,r$ occurs $r$ times.

Now consider the alphabet $A=\{a_{ij}\mid 1\leq i,j\leq r\}$ of
cardinality $r^2$. We convert the matrix $M$ to the matrix $M_{A}$
(shown in the right hand part of Fig.~\ref{matrices}) by replacing numbers
by letters according to the following rule: whenever the number
$i$ occurs in the column $j$ of $M$, we substitute it with the
letter $a_{ij}$ to get the corresponding entry in $M_A$.

Let $v_{t}$ be the word in the $t^{\mathrm{th}}$ row of the matrix
$M_A$. Consider the endomorphism $\gamma:A^+\rightarrow A^+$
defined by
$$\gamma(a_{ij})=v_{(i-1)r+j}.$$
Let $V_k$ be the set of all factors of the words in the sequence
$\{\gamma^m(a_{11})\}_{m=1,2,\dots}$ and let $0$ be a symbol
beyond $V_k$. We define a multiplication $\cdot$ on the set
$V_k\cup\{0\}$ as follows:
$$
u\cdot v=\left\{\begin{array}{cl}
uv&\ \text{if}\ u,v,uv\in V_k,\\
0 &\ \text{otherwise}.
\end{array}\right.
$$
Clearly, $\langle V_k\cup\{0\},\cdot\rangle$ becomes a semigroup
which we denote by $\mathcal{V}_k^0$. Using this semigroup, we can
conveniently reformulate two major combinatorial results by Sapir:

\begin{Prop} {\rm\cite[Proposition 2.1]{Sapir}}
\label{avoidable words} Let $X_k=\{x_1,\dots,x_k\}$ and $w\in
X_k^+$. Assume that there exists a homomorphism
$\varphi:X_k^+\to\mathcal{V}_k^0$ for which $\varphi(w)\ne0$. Then
there is an endomorphism $\psi:X_k^+\to X_k^+$ such that the word
$\psi(w)$ appears as a factor in the Zimin word $Z_k$.
\end{Prop}

\begin{Prop} {\rm\cite[Lemma 4.14]{Sapir}}
\label{avoidable identities} Let $X_k=\{x_1,\dots,x_k\}$ and
$w,w'\in X_k^+$. Assume that there exists a homomorphism
$\varphi:X_k^+\to\mathcal{V}_k^0$ for which
$\varphi(w)\ne\varphi(w')$. Then the identity $w=w'$ implies a
non-trivial semigroup identity of the form $Z_{k+1}=z$.
\end{Prop}

Now let $\overline{\mathcal{V}}_k^0$ denote the semigroup
anti-isomorphic to $\mathcal{V}_k^0$; we shall use the notation
$x\mapsto x^*$ for the mutual anti-isomorphisms between
$\mathcal{V}_k^0$ and $\overline{\mathcal{V}}_k^0$ in both
directions and denote $\{v^*\mid v\in V_k\}$ by $\overline{V}_k$.
Let
$$\mathcal{T}_k=\langle V_k\cup\overline{V}_k\cup\{0\},\cdot,{}^*\rangle$$
be the 0-direct union of $\mathcal{V}_k^0$ and
$\overline{\mathcal{V}}_k^0$; this means that we identify $0$ with
$0^*$, preserve the multiplication in both $\mathcal{V}_k^0$ and
$\overline{\mathcal{V}}_k^0$, and set $u\cdot v^*=u^*\cdot v=0$
for all $u,v\in V_k$. This is the unary semigroup we need.

It is clear that $\mathcal{T}_k$ is infinite and is generated (as
a unary semigroup) by the set $A$ which is finite. It remains to
verify that $\mathcal{T}_k$ satisfies every identity in at most
$k$ variables that holds in our initial unary semigroup
$\mathcal{S}$. So, let $p,q\in\FI(X_k)$ and suppose that
the identity $p=q$ holds in $\mathcal{S}$ but fails in
$\mathcal{T}_k$. Then there exists a unary semigroup homomorphism
$\varphi:\FI(X_k)\to\mathcal{T}_k$ for which
$\varphi(p)\ne\varphi(q)$. Hence, at least one of the elements
$\varphi(p)$ and $\varphi(q)$ is not equal to $0$; (without loss
of generality) assume that $\varphi(p)\ne 0$. Then we may also
assume $\varphi(p)\in V_k$; otherwise we may consider the identity
$p^*=q^*$ instead of $p=q$. Since $\varphi(p)\ne0$, there is no
letter $x\in X_k$ such that $p$ contains both $x$ and $x^*$. Now
we define a substitution
$\sigma:\FI(X_k)\to\FI(X_k)$ as follows:
$$\sigma(x)=\left\{\begin{array}{ll}
x^* & \text{if } p \text{ contains } x^*,\\
x & \text{otherwise}.\end{array}\right.$$ Then $\sigma(p)$ does
not contain any starred letter, thus being a plain word in
$X_k^+$. Since $\sigma^2$ is the identity mapping, we have
$\varphi(p)= (\varphi\sigma)(\sigma(p))$, and $\varphi\sigma$ maps
$X_k^+$ into $\mathcal{V}_k^0$. Now we consider two cases.

\emph{\textbf{Case 1:} $\sigma(q)$ contains a starred letter.} We
apply Proposition~\ref{avoidable words} to the plain word
$\sigma(p)$ and the semigroup homomorphism
$X_k^+\to\mathcal{V}_k^0$ obtained by restricting $\varphi\sigma$
to $X_k^+$. We conclude that there is an endomorphism $\psi$ of
$X_k^+$ such that the word $\psi(\sigma(p))$ appears as a factor
in the Zimin word $Z_k$. Thus, $Z_k=z'\psi(\sigma(p))z''$ for some
$z',z''$ (that may be empty). The endomorphism $\psi$ extends in a
natural way to an endomorphism of the free involutory semigroup
$\FI(X_k)$ and there is no harm in denoting the extension
by $\psi$ as well. The identity $p=q$ implies the identity
\begin{equation}
\label{consequence of p=q}
z'\psi(\sigma(p))z''=z'\psi(\sigma(q))z''.
\end{equation}
The left hand side of~\eqref{consequence of p=q} is $Z_k$ and the
identity is not trivial because its right hand side involves a
starred letter. Since $p=q$ holds in our initial semigroup
$\mathcal{S}$, so does~\eqref{consequence of p=q}. But this
contradicts the assumption that all Zimin words are involutory
isoterms for $\mathcal{S}$.

\emph{\textbf{Case 2:} $\sigma(q)$ contains no starred letter.} In
this case $\sigma(q)$ is a plain word in $X_k^+$, and we are in a
position to apply Proposition~\ref{avoidable identities} to the
semigroup identity $\sigma(p)=\sigma(q)$ and the semigroup
homomorphism $X_k^+\to\mathcal{V}_k^0$ obtained by restricting
$\varphi\sigma$ to $X_k^+$. We conclude that $\sigma(p)=\sigma(q)$
implies a non-trivial semigroup identity $Z_{k+1}=z$. Therefore
the identity $p=q$ implies $Z_{k+1}=z$, and we again get a
contradiction.
\end{proof}

Before passing to concrete examples of inherently \nfb\
involutory semigroups, we formulate a corollary of our proof
of Theorem~\ref{Theorem 2.2} which will be useful for our subsequent
paper(s). Proposition~\ref{avoidable words} easily implies that no word
in the sequence $\{\gamma^m(a_{11})\}_{m=1,2,\dots}$ has any square
(that is, a word of the form $ww$) as a factor. Hence, the semigroup
$\mathcal{V}_k^0$ satisfies the identity\footnote{Strictly speaking
the following expression is not an identity as defined before (since
0 is not a term) but rather an abbreviation for the identities
$x^2y=x^2=yx^2$. However, referring to such abbreviations as identities
is a standard convention which we adopt.} $x^2=0$. This identity is
clearly inherited by the involutory semigroup $\mathcal{T}_k$ which
by its construction satisfies also the identity $xx^*=0$. Since $\mathcal{T}_k$
is finitely generated and infinite, we arrive at the following conclusion:
\begin{Cor}
\label{new corollary}
If a variety $\mathbf{V}$ of involutory semigroups satisfies no
non-trivial identity of the form $Z_{k+1}=z$ and all members of
$\mathbf{V}$ satisfying the identities $xx^*=x^2=0$ are locally
finite, then $\mathbf{V}$ is not finitely based.
\end{Cor}
This result is parallel to \cite[Proposition~3]{sapirburnside} in
the plain semigroup case.

Now consider the \emph{twisted Brandt monoid}
$\TB=\langle B_2^1,\cdot,{}^*\rangle$, where $B_2^1$ is the set
of the following six $2\times 2$-matrices:
$$\begin{pmatrix} 0 & 0\\ 0 & 0\end{pmatrix},\
\begin{pmatrix} 1 & 0\\ 0 & 0\end{pmatrix},\
\begin{pmatrix} 0 & 1\\ 0 & 0\end{pmatrix},\
\begin{pmatrix} 0 & 0\\ 1 & 0\end{pmatrix},\
\begin{pmatrix} 0 & 0\\ 0 & 1\end{pmatrix},\
\begin{pmatrix} 1 & 0\\ 0 & 1\end{pmatrix},$$
the binary operation $\cdot$ is the usual matrix multiplication
and the unary operation ${}^*$ fixes the matrices
$$\begin{pmatrix} 0 & 0\\ 0 & 0\end{pmatrix},\
\begin{pmatrix} 0 & 1\\ 0 & 0\end{pmatrix},\
\begin{pmatrix} 0 & 0\\ 1 & 0\end{pmatrix},\
\begin{pmatrix} 1 & 0\\ 0 & 1\end{pmatrix}$$
and swaps each of the matrices
$$\begin{pmatrix} 1 & 0\\ 0 & 0\end{pmatrix},\
\begin{pmatrix} 0 & 0\\ 0 & 1\end{pmatrix}$$
with the other one.

\begin{Cor}
\label{twisted Brandt} The twisted Brandt monoid $\TB$ is
inherently \nfb.
\end{Cor}

\begin{proof}
By Theorem~\ref{Theorem 2.2} we only have to show that $\TB$
satisfies no non-trivial involutory semigroup identity of the form
$Z_n=z$. If $z$ is a plain semigroup word, we can refer to
\cite[Lemma~3.7]{sapirburnside} which shows that the semigroup
$\langle B_2^1,\cdot\rangle$ does not satisfy any non-trivial
\textbf{semigroup} identity of the form $Z_n=z$. If we suppose
that the involutory word $z$ contains a starred letter, we can
substitute the matrix $\left(\begin{smallmatrix}1 & 0\\ 0 &
0\end{smallmatrix}\right)$ for all letters occurring in $Z_n$ and
$z$. Since this matrix is idempotent, the value of the word $Z_n$
under this substitution equals $\left(\begin{smallmatrix}1 & 0\\ 0
& 0\end{smallmatrix}\right)$. On the other hand, $z$ evaluates to
a product involving the matrix $\left(\begin{smallmatrix}1 & 0\\ 0
& 0\end{smallmatrix}\right)^*= \left(\begin{smallmatrix}0 & 0\\ 0
& 1\end{smallmatrix}\right)$, and it is easy to see that such a
product is equal to either $\left(\begin{smallmatrix}0 & 0\\ 0 &
1\end{smallmatrix}\right)$ or $\left(\begin{smallmatrix}0 & 0\\ 0
& 0\end{smallmatrix}\right)$. Thus, the identity $Z_n=z$ cannot
hold in $\TB$ in this case as well.
\end{proof}

An equivalent way to define $\TB$ is to consider the 5-element
unary \Rm\ semigroup over the trivial group $\mathcal{E}=\{e\}$
with the sandwich matrix
$$\begin{pmatrix}
0 & e\\
e & 0
\end{pmatrix}$$
and then to adjoin to this unary \Rm\ semigroup an identity
element. For convenience and later use we note that $\TB$ can thus
be realized as the set
$$\{(1,1),(1,2),(2,1),(2,2),0,1\}$$ endowed with the operations
\begin{gather}
\label{operations in TB}
(i,j)\cdot(k,\ell)=\left\{\begin{array}{cl}
(i,\ell)&\ \text{if}\ (j,k)\in \{(1,2),(2,1)\},\\
0 &\ \text{otherwise};
\end{array}\right.\\
1\cdot x=x=x\cdot 1,\ 0\cdot x=0=x\cdot 0 \text{ for all }x; \notag\\
 (i,j)^* = (j,i),\ 1^*=1,\ 0^*=0\notag.
\end{gather}

Suppose that $\mathcal{S}$ is a finite unary semigroup for
which the variety $\var\mathcal{S}$ contains an inherently
\nfb\ unary semigroup. Immediately from the definition it
follows that $\mathcal{S}$ is also inherently \nfb. This
observation is useful, in particular, for the justification
of our second example of an involutory inherently \nfb\ semigroup.
This is a `twisted version' $\TA$ of another 6-element semigroup
that often shows up under the name $A_2^1$ in the theory of semigroup
varieties. The unary semigroup $\TA$ is formed by the 6 matrices
$$\begin{pmatrix} 0 & 0\\ 0 & 0\end{pmatrix},\
\begin{pmatrix} 1 & 0\\ 0 & 0\end{pmatrix},\
\begin{pmatrix} 0 & 1\\ 0 & 0\end{pmatrix},\
\begin{pmatrix} 1 & 0\\ 1 & 0\end{pmatrix},\
\begin{pmatrix} 0 & 1\\ 0 & 1\end{pmatrix},\
\begin{pmatrix} 1 & 0\\ 0 & 1\end{pmatrix}$$
under the usual matrix multiplication and the unary operation that
swaps each of the matrices
$$\begin{pmatrix} 1 & 0\\ 0 & 0\end{pmatrix},\
\begin{pmatrix} 0 & 1\\ 0 & 1\end{pmatrix}$$
with the other one and fixes all other matrices. Alternatively,
$\TA$ is obtained from the 5-element unary \Rm\ semigroup
$\mathcal{A}_2$ over $\mathcal{E}=\{e\}$ with the sandwich matrix
\begin{equation}
\label{matrix for TA}
\begin{pmatrix}
0 & e\\
e & e
\end{pmatrix}
\end{equation}
by adjoining an identity element. Again, for later use, we note
that $\TA$ can be realized as the set
$$\{(1,1),(1,2),(2,1),(2,2),0,1\}$$ endowed with the operations
\begin{gather}
\label{operations in TA}
(i,j)\cdot(k,\ell)=\left\{\begin{array}{cl}
(i,\ell)&\ \text{if}\ (j,k)\ne (1,1)\\
0 &\ \text{if}\ (j,k)=(1,1);
\end{array}\right.\\
1\cdot x=x=x\cdot 1,\ 0\cdot x=0=x\cdot 0 \text{ for all }x; \notag\\
 (i,j)^* = (j,i),\ 1^*=1,\ 0^*=0\notag.
\end{gather}

\begin{Cor}
\label{twisted A} The involutory semigroup $\TA$ is inherently
\nfb.
\end{Cor}

\begin{proof}
We represent $\TA$ as in \eqref{operations in TA} and $\TB$ as in
\eqref{operations in TB} and consider the direct square
$\TA\times\TA$. It is then easy to check that the twisted Brandt
monoid $\TB$ is a homomorphic image of the unary subsemigroup of
$\TA\times\TA$ generated by the pairs $(1,1)$,
$\bigl((1,1),(2,2)\bigr)$ and $\bigl((2,2),(1,1)\bigr)$. Thus,
$\TB$ belongs to $\var\TA$. Since by Corollary~\ref{twisted
Brandt} $\TB$ is inherently \nfb, so is $\TA$.
\end{proof}

Sapir~\cite[Proposition~7]{sapirburnside} has shown that a (plain)
finite semigroup $\mathcal{S}$ is inherently \nfb\ \textbf{if
and only if} all Zimin words are isoterms for $\mathcal{S}$, that
is, $\mathcal{S}$ satisfies no non-trivial semigroup identity of
the form $Z_n=z$. Our Theorem~\ref{Theorem 2.2} models the `if'
part of this statement but we do not know whether or not the `only
if' part transfers to the involutory environment. Some partial
results in this direction have been recently obtained by the
second author \cite{Dolinka}. Here we present yet another special
result which however suffices for our purposes.

\begin{Prop}
\label{NINFB} Let $\mathcal{S}=\langle S,\cdot,{}^*\rangle$ be a
finite involutory semigroup and suppose that there exists an
involutory word $\om(x)$ in one variable $x$ such that
$\mathcal{S}$ satisfies the identity $x=x\om(x)x$. Then
$\mathcal{S}$ is not inherently \nfb.
\end{Prop}

\begin{proof}
We have to construct a finite set of identities that defines
a locally finite variety of involutory semigroups containing
$\mathcal{S}$. The crucial step towards this goal consists in
verifying that the identity $x=x\om(x)x$ allows one to express
right divisibility in terms of equational logic. This being done,
we shall be in a position to closely follow powerful arguments
developed by Margolis and Sapir in~\cite{margolissapir}.

We say that elements $a,b\in S$ \emph{divide each other on the right}
and write $a\mathrel{\Rc}b$ if either $a=b$ or there exist $s,t\in S$
such that $a=bs$ and $b=at$. We say that $b$ \emph{strictly divides}
$a$ and write $a\mathrel{{<}_{\Rc}}b$ if $a=bs$ for some $s\in S$
but $b\ne a$ and $b\ne at$ for any $t\in S$. Clearly, $\Rc$ is an
equivalence relation (known as the \emph{right Green relation} in
semigroup theory) and $<_{\Rc}$ is transitive and anti-reflexive.

Since $\mathcal{S}$ satisfies the identity $x=x\om(x)x$, we have
$a\mathrel{\Rc}a\om(a)$ for each element $a\in S$. (Indeed,
$a=a\om(a)\cdot a$ and $a\om(a)=a\cdot\om(a)$.) Thus, for
$a,b\in S$, we have $a\mathrel{\Rc}b$ if and only if
$a\om(a)\mathrel{\Rc}b\om(b)$. Since $a\om(a)$ and $b\om(b)$ are
idempotents, that latter condition is equivalent to the two
equalities $a\om(a)\cdot b\om(b)=b\om(b)$ and $b\om(b)\cdot
a\om(a)=a\om(a)$. In particular, for $u,v\in S$ we have
$uv\mathrel{\Rc}u$ if and only if $uv\om(uv)\cdot u\om(u)=u\om(u)$
(since the second equality $u\om(u)\cdot uv\om(uv)=uv\om(uv)$ is
always true).

Let $Z_n'$ be the word obtained from the Zimin word $Z_n$ by
deleting the last letter (which is $x_1$), that is, $Z_n'x_1=Z_n$.
Further let $h$ denote the length of the longest possible chain
of the form
$$s_1\mathrel{{<}_{\Rc}}s_2\mathrel{{<}_{\Rc}}\cdots\mathrel{{<}_{\Rc}}s_k.$$
Set $n=h+1$; Lemma~7 in~\cite{margolissapir} shows that under every
evaluation of the letters $x_1,\dots,x_n$ in $\mathcal{S}$, the values
of the words $Z_n'$ and $Z_n$ divide each other on the right. As explained
in the previous paragraph, this can be restated as saying that the values
of the terms $Z_n\om(Z_n)\cdot Z_n'\om(Z_n')$ and $Z_n'\om(Z_n')$
coincide under every evaluation of $x_1,\dots,x_n$ in $\mathcal{S}$, that is,
$\mathcal{S}$ satisfies the identity
\begin{equation}
\label{ziminR} Z_n\om(Z_n)\cdot Z_n'\om(Z_n')=Z_n'\om(Z_n').
\end{equation}
On the other hand, in each involutory semigroup $\mathcal{T}=\langle T,\cdot,{}^*\rangle$
which satisfies $x=x\om(x)x$ and \eqref{ziminR}, the values of the words $Z_n'$
and $Z_n$ under every evaluation of $x_1,\dots,x_n$ necessarily divide
each other on the right. This implies that such $\mathcal{T}$ satisfies
the implication
\begin{equation}
\label{quasidentity}
xZ_n=yZ_n\rightarrow xZ_n'=yZ_n'.
\end{equation}
Indeed, suppose that under some evaluation $\varphi$ of the letters
$x,y,x_1,\dots,x_n$ in $\mathcal{T}$, the words $xZ_n$ and $yZ_n$
happen to take the same value, that is, $\varphi(xZ_n)=\varphi(yZ_n)$.
Since $\varphi(Z_n)\mathrel{\Rc}\varphi(Z_n')$, there exists $t\in T$
such that $\varphi(Z_n)t=\varphi(Z_n')$. Hence
$$\varphi(xZ_n')=\varphi(xZ_n)t=\varphi(yZ_n)t=\varphi(yZ_n'),$$
that is, the words $xZ_n'$ and $yZ_n'$ also take a common value under $\varphi$.

Lemma 8 in \cite{margolissapir} shows that a finitely generated semigroup
satisfying \eqref{quasidentity} is finite if and only if it satisfies
the identity
\begin{equation}
\label{periodic}
x^k=x^{k+\ell}
\end{equation}
for some $k,\ell\ge 1$ and has only locally finite subgroups. We note that
an involutory semigroup $\mathcal{T}=\langle T,\cdot,{}^*\rangle$ is finitely
generated if and only if so is the semigroup $\langle T,\cdot\rangle$.
An identity of the form \eqref{periodic} definitely holds in $\mathcal{S}$
since $\mathcal{S}$ is finite. Hence it suffices to find a finite number of
identities which hold in $\mathcal{S}$ and which force each (involutory) semigroup
to have only locally finite subgroups.

We can proceed as at the end of \cite{margolissapir}. Let
$\mathcal{G}$ be the direct product of all maximal subgroups of
$\mathcal{S}$. By the Oates-Powell theorem~\cite{oatespowell}, see
also \cite[\S5.2]{Ne}, the locally finite variety
$\var\mathcal{G}$ generated by the finite group $\mathcal{G}$ can
be defined by a single identity $v(x_1,\dots,x_m)=1$. The left
hand side $v$ of this identity can be assumed to contain no occurrence
of the inversion ${}^{-1}$, that is, we may assume that $v$ is a plain semigroup
word in the letters $x_1,\dots,x_m$. Now let $\mathcal{F}=\mathcal{F}(x_1,\dots,x_m)$
be the $m$-generated relatively free semigroup in the (locally finite)
semigroup variety generated by the semigroup $\langle S,\cdot\rangle$.
The semigroup $\mathcal{F}$ is finite and hence has a least
ideal; this ideal is known to be a union of (isomorphic)
subgroups (Sushkevich's theorem, see~\cite{CP}). Let
$\mathcal{H}$ be any of these subgroups and let $e$ be the identity
element of $\mathcal{H}$. We denote by $u(x_1,\dots,x_m)$ a word
whose value in $\mathcal{F}$ is $e$. Since $e^2=e$, we see that
$\mathcal{F}$ (and therefore $\mathcal{S}$) satisfies the identity
\begin{equation}
\label{idempotent law} u=u^2.
\end{equation}
For every element $g\in\mathcal{F}$, the product $ege$ belongs to
$\mathcal{H}$. As observed in \cite{KR}, the group $\mathcal{H}$
belongs to the variety $\var\mathcal{G}$. Consequently, $\mathcal{F}$
(and therefore $\mathcal{S}$) satisfies the identity
\begin{equation}
\label{crucial law} v(ux_1u,\dots,ux_mu)=u.
\end{equation}
Note that both sides of \eqref{crucial law} are plain semigroup words in
the letters $x_1,\dots,x_m$.

Now consider the variety $\mathbf{V}$ of involutory semigroups defined
by the identity $x=x\om(x)x$, an identity of the form \eqref{periodic}
holding in $\mathcal{S}$, and the identities \eqref{ziminR},
\eqref{idempotent law}, and \eqref{crucial law}. Since by the
construction $\mathcal{S}$ satisfies all the listed identities,
$\mathcal{S}$ is a member of $\mathbf{V}$. Let
$\mathcal{T}=\langle T,\cdot,{}^*\rangle$ be any finitely generated member of $\mathbf{V}$;
then, as already mentioned, the semigroup $\langle T,\cdot\rangle$ is also
finitely generated. The first and the third identity ensure that the
semigroup $\langle T,\cdot\rangle$ satisfies the implication \eqref{quasidentity},
and therefore, it is finite provided that all its subgroups are locally finite.
Each group that satisfies the identities \eqref{idempotent law}
and \eqref{crucial law} satisfies the identity $v(x_1,\dots,x_m)=1$, whence
this group belongs to $\var\mathcal{G}$ and so is locally finite. Altogether,
$\mathcal{T}$ is finite. Thus, $\mathbf{V}$ is locally finite and \fb, and
the proposition is proved.
\end{proof}

Proposition \ref{NINFB} implies in particular that no finite regular $*$-semigroup can be inherently \nfb\ as one
can use $x^*$ in the role of the term $\om(x)$. In particular, the unary semigroup $\langle B_2^1,\cdot,{}^T\rangle$,
where the unary operation is the usual matrix transposition, is not inherently \nfb\ (this fact was first discovered
by Sapir, see~\cite{sapirinverse}), even though it is not finitely based~\cite{kleiman}.

\section{Applications}\label{applications}

\subsection{Matrix semigroups with Moore-Penrose inverse}
\label{Moore-Penrose}

Certainly, the most common unary operation for matrices is
transposition. However, it is convenient for us to start with
analyzing matrix semigroups with Moore-Penrose inverse because
this analysis will help us in considering semigroups with
transposition.

We first recall the notion of Moore-Penrose inverse. This has been
discovered by Moore \cite{moore} and independently by Penrose
\cite{P} for complex matrices, but has turned out to be a fruitful
concept in a more general setting---see \cite{BIG} for a
comprehensive treatment.

The following results were obtained by Drazin \cite{D}.

\begin{Prop}\label{Theorem 3.3} {\rm\cite[Proposition 1]{D}}
Let $\mathcal{S}$ be an involutory semigroup. Then, for any given
$a\in\mathcal{S}$, the four equations
\begin{equation}
\label{penrose} axa=a,\ xax=x,\ (ax)^*=ax,\ (xa)^*=xa
\end{equation}
have at most one common solution $x\in\mathcal{S}$.
\end{Prop}

For an element $a$ of  an involutory semigroup $\mathcal{S}$, we
denote by $a^\dag$ the unique common solution $x$ of the equations
\eqref{penrose}, provided it exists, and call $a^\dag$ the
\emph{Moore-Penrose inverse} of $a$.

Recall that an element $a\in\mathcal{S}$ is said to be
\emph{regular}, if there is an $x\in\mathcal{S}$ such that
$axa=a$. Concerning existence of the Moore-Penrose inverse, we
have the following
\begin{Prop}\label{Theorem 3.4} {\rm\cite[Proposition 2]{D}}
Let $\mathcal{S}$ be an involutory semigroup satisfying the
implication
\begin{equation}
\label{3.1} x^*x=x^*y=y^*x=y^*y\rightarrow x=y.
\end{equation}
Then for an arbitrary $a\in\mathcal{S}$, the Moore-Penrose inverse
$a^\dag$ exists if and only if $a^*a$ and $aa^*$ are regular
elements.
\end{Prop}

Let $\langle R,+,\cdot\rangle$ be a ring. An \emph{involution of
the ring} is an involution $x\mapsto x^*$ of the semigroup
$\langle R,\cdot\rangle$ satisfying in addition the identity
$(x+y)^*=x^* +y^*$. For ring involutions, the
implication~\eqref{3.1} is easily seen to be equivalent to
\begin{equation}
\label{3.2} x^*x=0\rightarrow x=0.
\end{equation}
Now suppose that $\mathcal{K}=\langle K,+,\cdot\rangle$ is a field
that admits an involution $x\mapsto\ol x$. Then the matrix ring
$\mathrm{M}_n(\mathcal{K})$ has an involution that naturally
arises from the involution of $\mathcal{K}$, namely
$(a_{ij})\mapsto(a_{ij})^*:=(\ol{a_{ij}})^T$. This involution of
$\mathrm{M}_n(\mathcal{K})$ in general does not satisfy the
implication~\eqref{3.2}. However, it does satisfy \eqref{3.2}
if and only if the equation
\begin{equation} \label{3.3}
x_1\ol{x_1}+x_2\ol{x_2}+\dots+x_n\ol{x_n}=0
\end{equation}
admits only the trivial solution $(x_1,\dots,x_n)=(0,\dots,0)$ in
$K^n$. Since all elements of $\mathrm{M}_n(\mathcal{K})$ are
regular, this means that the Moore-Penrose inverse exists ---
subject to the  involution $(a_{ij})\mapsto
(a_{ij})^*=(\ol{a_{ij}})^T$ --- whenever \eqref{3.3} admits only
the trivial solution. (The classical Moore-Penrose inverse is
thereby obtained by putting $\mathcal{K}=\bb C$, the field of
complex numbers, endowed with the usual complex conjugation
$z\mapsto \ol z$.) On the other hand, it is easy to see that the
condition that \eqref{3.3} has only the trivial solution is
necessary: if $(a_1,\dots,a_n)$ were a non-trivial solution to
\eqref{3.3}, then the matrix formed by $n$ identical rows
$(a_1,\dots,a_n)$ would have no Moore-Penrose inverse.

The proof of the main result of this subsection requires an
explicit calculation of the Moore-Penrose inverses of certain
rank~1 matrices. Thus, we present a simple method for such a
calculation. For a row vector $a=(a_1,\dots,a_n)\in K^n$, where
$\mathcal{K}=\langle K,+,\cdot\rangle$ is a field with an
involution $x\mapsto\ol x$, let $a^*$ denote the column vector
$(\ol{a_1},\dots,\ol{a_n})^T$. It is easy to see that any $n\times
n$-matrix $A$ of rank~1 over $\mathcal{K}$ can be represented as
$A=b^*c$ for some non-zero row vectors $b,c\in K^n$. Provided that
\eqref{3.3} admits only the trivial solution in $K^n$, one gets
$A^\dag$ as follows:
\begin{equation}
\label{calculating MP inverse for rank 1}
A^\dag=c^*(cc^*)^{-1}(bb^*)^{-1}b.
\end{equation}
Here $bb^*$ and $cc^*$ are non-zero elements of $\mathcal{K}$
whence their inverses in $\mathcal{K}$ exist. In order to
justify~\eqref{calculating MP inverse for rank 1}, it suffices to
check that the right hand side of~\eqref{calculating MP inverse
for rank 1} satisfies the simultaneous equations~\eqref{penrose}
with the matrix $A$ in the role of $a$, and this is
straightforward. Note that formula \eqref{calculating MP inverse
for rank 1} immediately shows that $A^\dag$ is a scalar multiple
of $A^*=c^*b$, namely
\begin{equation}\label{MP is scalar multiple}
A^\dag=\frac{1}{cc^*\cdot bb^*}A^*.
\end{equation}

So, we can formulate one of the highlights of the section --- a
result that reveals an unexpected feature of a rather classical
and well studied object.
\begin{Thm} \label{Theorem 3.5}
Let $\mathcal{K}=\langle K,+,\cdot\rangle$ be a field having an
involution $x\mapsto \ol x$ for which the equation $x\ol x +y\ol y
=0$ has only the trivial solution $(x,y)=(0,0)$ in $K^2$. Then the
unary semigroup
$\langle\mathrm{M}_2(\mathcal{K}),\cdot,{}^\dag\rangle$ of all
$2\times 2$-matrices over $\mathcal{K}$ endowed with Moore-Penrose
inversion $^\dag$ --- subject to the involution $(a_{ij})\mapsto
(a_{ij})^*= (\ol{a_{ij}})^T$
--- has no finite basis of identities.
\end{Thm}

\begin{proof}
Set $\Sc=\langle\mathrm{M}_2(\mathcal{K}),\cdot,{}^\dag\rangle$.
By Theorem~\ref{Theorem 2.1} and Lemma~\ref{Lemma 3.1} it is
sufficient to show that
\renewcommand{\labelenumi}{\theenumi)}
\begin{enumerate}
\item $\mathcal{K}_3\in\var\Sc$,
\item there exists a group $\mathcal{G}\in\var\Sc$
such that $\mathcal{G}\notin\var\H(\Sc)$.
\end{enumerate}

In order to prove 1), consider the following sets of rank~1
matrices in $\mathrm{M}_2(\mathcal{K})$:
\begin{eqnarray}
\notag H_{11}=\left\{\begin{pmatrix} x & x\\ x &
x\end{pmatrix}\right\}, & H_{12}=\left\{\begin{pmatrix} x & 0\\ x
& 0\end{pmatrix}\right\}, & H_{13}=\left\{\begin{pmatrix}
0 & x\\ 0 & x\end{pmatrix}\right\}, \\
\label{C3 is in M2} H_{21}=\left\{\begin{pmatrix} x & x\\ 0 &
0\end{pmatrix}\right\}, & H_{22}=\left\{\begin{pmatrix} x & 0\\ 0
& 0\end{pmatrix}\right\}, & H_{23}=\left\{\begin{pmatrix}
0 & x\\ 0 & 0\end{pmatrix}\right\},\\
\notag H_{31}=\left\{\begin{pmatrix} 0 & 0\\ x &
x\end{pmatrix}\right\}, & H_{32}=\left\{\begin{pmatrix} 0 & 0\\ x
& 0\end{pmatrix}\right\}, & H_{33}=\left\{\begin{pmatrix} 0 & 0\\
0 & x\end{pmatrix}\right\},
\end{eqnarray}
where in each case $x$ runs over $K\setminus\{0\}$. Observe that
$\mathcal{K}$ cannot be of characteristic $2$, since the equation
$x\ol x+y\ol y=0$ has only the trivial solution in ${K}^2$. Taking
this into account, a straightforward calculation shows that
\begin{equation}
\label{multiplication in M2} H_{ij}\cdot
H_{k\ell}=\left\{\begin{array}{cl}
H_{i\ell}&\ \text{if}\ (j,k)\ne(2,3),(3,2),\\
0 &\ \text{otherwise}.
\end{array}\right.
\end{equation}
Hence the set
$$T=\bigcup_{1\le i,j\le3}H_{ij}\cup\{0\}$$
is closed under multiplication so that this set forms a
subsemigroup $\mathcal{T}$ of $\Sc$ and the partition $\Hc$ of $T$
into the classes $H_{ij}$ and $\{0\}$ is a congruence on
$\mathcal{T}$. Equation \eqref{MP is scalar multiple} shows that
\begin{equation}
\label{inversion in M2} H_{ij}^\dag=H_{ji}.
\end{equation}
We see that $T$ is closed under Moore-Penrose inversion and $\Hc$
respects $^\dag$, thus is a congruence on the unary semigroup
$\mathcal{T}'=\left<T,\cdot,{}^\dag\right>$. Now comparing
\eqref{multiplication in M2} and \eqref{inversion in M2} with the
multiplication and inversion rules in $\mathcal{K}_3$ (see
\eqref{operations in C3}), we conclude that $\mathcal{T}'/\Hc$ and
$\mathcal{K}_3$ are isomorphic as unary semigroups. Hence
$\mathcal{K}_3$ is in $\var\Sc$.

For 2) we merely let $\mathrm{GL}_2(\mathcal{K})$, the group of
all invertible $2\times 2$-matrices over $\mathcal{K}$, play the
role of $\mathcal{G}$. Since Moore-Penrose inversion on
$\mathrm{GL}_2(\mathcal{K})$ coincides with usual matrix
inversion, we observe that $\mathrm{GL}_2(\mathcal{K})$  is a
unary subsemigroup of $\Sc$. Moreover, since $AA^\dag$ is the
identity matrix for every invertible matrix $A$, we conclude that,
with the exception of the identity matrix, the Hermitian
subsemigroup $\H(\Sc)$ contains only matrices of rank~1, that is,
$\H(\Sc)\subseteq\mathrm{L}^1_2(\mathcal{K})$, the set
of all matrices of rank at most~1 with the identity matrix
adjoined. By Corollary~\ref{identity for rank 1} the semigroup
$\langle\mathrm{L}^1_2(\mathcal{K}),\cdot\rangle$ satisfies the
identity $x^2yx=xyx^2$. Consequently, each group in $\var\H(\Sc)$ is
abelian, while the group $\mathrm{GL}_2(\mathcal{K})$ is
non-abelian. Thus, $\mathrm{GL}_2(\mathcal{K})$ is contained in
$\var\Sc$ but is not contained in $\var\H(\Sc)$, as required.
\end{proof}

\begin{Rmk}\label{remark on MP inverse}
Apart from any subfield of $\bb C$ closed under complex
conjugation, Theorem~\ref{Theorem 3.5} applies, for instance, to
finite fields $\mathcal{K}=\langle K,+,\cdot\rangle$ for which
$|K|\equiv 3\!\pmod{4}$, endowed with the trivial involution
$x\mapsto \ol x=x$; the latter follows from the fact that the
equation $x^2+1=0$ admits no solution in $\mathcal{K}$ if and only
if $|K|\equiv 3\!\pmod{4}$ (cf.\ \cite[Theorem
3.75]{LidlNiederreiter}). Moreover, by slightly changing the
arguments one can show an analogous result for $\mathcal{K}$ being
any skew-field of quaternions closed under conjugation.
\end{Rmk}

The reader may ask whether or not the restriction on the size of
matrices is essential in Theorem~\ref{Theorem 3.5}. For some
fields, it definitely is. For instance, for finite fields with the
trivial involution $x\mapsto \ol x=x$, no extension of
Theorem~\ref{Theorem 3.5} to $n\times n$-matrices with $n>2$ is
possible simply because the Moore-Penrose inverse is only a
partial operation in this case. Indeed, it is a well known
corollary of the Chevalley-Warning theorem (see, e.g.,
\cite[Corollary 2 in \S1.2]{Serre}) that the equation
$x_1^2+\dots+x_n^2=0$ (that is \eqref{3.3} with the trivial
involution) admits a non-trivial solution in any finite field
whenever $n>2$.

The situation is somewhat more complicated for subfields of $\bb
C$. Theorem~\ref{Theorem 2.1} does not apply here because of the
following obstacle. It is well known (see, for example,
\cite[p.\,101]{mks}) that the two matrices
\begin{equation}
\label{free subgroup} \zeta=\begin{pmatrix} 1 & 0\\ 2 &
1\end{pmatrix} \ \text{ and } \ \eta=\begin{pmatrix} 1 & 2\\ 0 &
1\end{pmatrix}
\end{equation}
generate a free subgroup of $\left<\mathrm{SL}_2(\bb
Z),\cdot,{}^{-1}\right>$. On the other hand, it is easy to verify
that for any subfield $\mathcal{K}$ of $\bb C$ closed under
complex conjugation, the mapping $\langle\mathrm{SL}_2(\bb
Z),\cdot,{}^{-1}\rangle\to\langle\mathrm{M}_3(\mathcal{K}),\cdot,{}^\dag\rangle$
defined by $A\mapsto\left(\begin{matrix} \raisebox{-4.5pt}{$A$}\\
\begin{smallmatrix} 0 & 0 \end{smallmatrix}\end{matrix}
\begin{smallmatrix} 0 \\ 0\\ 0\end{smallmatrix}\right)$
is an embedding of unary semigroups. Thus, for $n>2$, the unary
semigroup
$\langle\mathrm{Sing}_n(\mathcal{K}),\cdot,{}^\dag\rangle$ of
singular $n\times n$-matrices contains a free non-abelian group,
whence every group belongs to the unary semigroup variety
generated by $\langle\mathrm{Sing}_n(\mathcal{K}),\cdot,{}^\dag\rangle$.
Now we observe that $\mathrm{Sing}_n(\mathcal{K})$ is contained in
the Hermitian subsemigroup of $\langle\mathrm{M}_n(\mathcal{K}),\cdot,{}^\dag\rangle$.
Indeed, it was proved in \cite{Erdos} (see also \cite{AM} for a recent
elementary proof) that the semigroup
$\langle\mathrm{Sing}_n(\mathcal{K}),\cdot\rangle$ is generated
by idempotent matrices. For an arbitrary idempotent matrix
$A\in\mathrm{M}_n(\mathcal{K})$, let $$N(A)=\{x\in K^n\mid xA=0\}\
\text{ and } \  F(A)=\{x\in K^n\mid xA=x\}$$ be the null-space and
the fixed-point-space of $A$, respectively. Now consider two
matrices of orthogonal projectors: $P_1$, the matrix of the
orthogonal projector to the space $F(A)$, and $P_2$, the matrix of
the orthogonal projector to the space $N(A)^\perp$. As any
orthogonal projector matrix $P$ satisfies $P=P^2=P^\dag$, both
$P_1=P_1P_1^\dag$ and $P_2=P_2P_2^\dag$ belong to the Hermitian
subsemigroup $\H(\mathrm{M}_n(\mathcal{K}))$, but then $A$ also
belongs to $\H(\mathrm{M}_n(\mathcal{K}))$ since
$A=(P_1P_2)^\dag$, see \cite[Exercise~5.15.9a]{Meyer}. Thus,
$\mathrm{Sing}_n(\mathcal{K})\subseteq
\H(\mathrm{M}_n(\mathcal{K}))$, whence no group $\mathcal{G}$ can
satisfy the condition of Theorem~\ref{Theorem 2.1}.

However, the fact that Theorem~\ref{Theorem 2.1} cannot be applied
to, say, the unary semigroup $\langle\mathrm{M}_3(\bb
C),\cdot,{}^\dag\rangle$ does not yet mean that the identities of
this semigroup are finitely based. We thus have the following open
question.
\begin{Problem}
\label{problem on Moore-Penrose} Is the unary semigroup
$\langle\mathrm{M}_n(\mathcal{K}),\cdot,{}^\dag\rangle$ not
finitely based for each subfield $\mathcal{K}$ of $\bb C$ closed
under complex conjugation and for all $n>2$?
\end{Problem}

It is known that $\langle\mathrm{M}_n(\mathcal{K}),\cdot,{}^\dag\rangle$
satisfies rather involved identities (see \cite{Cline} for an example),
so the conjecture that these identities admit no finite basis looks
quite natural. In connection with Problem \ref{problem on Moore-Penrose},
we also observe that the proofs of Theorem~\ref{Theorem 3.5} and
Corollary~\ref{identity for rank 1} readily yield the following:
\begin{Rmk}
\label{rank 1 + invertible} For each conjugation-closed subfield
$\mathcal{K}$ of $\bb C$  and for all $n>2$,  the unary semigroup
$\langle\mathrm{L}_n(\mathcal{K})\cup G,\cdot,{}^\dag\rangle$ consisting of all
matrices of rank at most~1 and all matrices from some non-abelian subgroup
$\langle G,\cdot,{}^{-1}\rangle$ of
$\langle\mathrm{GL}_n(\mathcal{K}),\cdot,{}^{-1}\rangle$ has no finite identity
basis.
\end{Rmk}

Another natural related structure is the semigroup
$\textrm{M}_n(\mathcal{K})$ endowed with \textbf{both} unary
operations ${}^\dag$ and ${}^*$. Recall that the Moore-Penrose
inverse is in fact defined in terms of identities involving
both these operations---namely, Proposition~\ref{Theorem 3.4}
implies that $\langle\mathrm{M}_n(\mathcal{K}),\cdot,{}^\dag,{}^*\rangle$
satisfies the identities
\begin{equation}
\label{penrose1} xx^\dag x=x,\ x^\dag xx^\dag=x^\dag,\ (xx^\dag)^*=xx^\dag,\ (x^\dag x)^*=x^\dag x
\end{equation}
and these identities uniquely determine the operation $A\mapsto A^\dag$.
This might have provoked one to conjecture that the identities \eqref{penrose1}
together with the identities $(xy)^* = y^*x^*$ and $(x^*)^* = x$ form
a basis for $\Id\langle\mathrm{M}_n(\mathcal{K}),\cdot,{}^\dag,{}^*\rangle$.
However our techniques show that this is not the case at least for $n=2$.
\begin{Thm}\label{both operations}
Let $\mathcal{K}$ be a field as in Theorem~$\ref{Theorem 3.5}$;
then $\langle\mathrm{M}_2(\mathcal{K}),\cdot,{}^\dag,{}^*\rangle$
is not finitely based as an algebraic structure of type $(2,1,1)$.
\end{Thm}

\begin{proof}
The characteristic of $\mathcal{K}$ is not $2$ whence the group
$$\mathcal{G}=\{A\in \textrm{GL}_2(\mathcal{K})\mid A^\dag=A^*\}$$
is non-abelian. Indeed, on the prime subfield of $\mathcal{K}$,
the involution $x \mapsto \ol x$ is the identity automorphism; so,
for matrices over the prime subfield, conjugation ${}^*$ coincides
with transposition,
and thus, for example, $(\begin{smallmatrix} 0 & 1\\
1 & 0\end{smallmatrix})$ and $(\begin{smallmatrix} 0 & -1\\ 1 & 0
\end{smallmatrix})$ are two non-commuting members of $\mathcal{G}$. Set
$\Ac=\langle\mathrm{M}_2(\mathcal{K}),\cdot,{}^\dag,{}^*\rangle$;
the algebraic structure $\langle G,\cdot,{}^\dag,{}^*\rangle$,
that is, the group $\mathcal{G}$ with inversion taken twice as
unary operation, belongs to $\var\Ac$. Now, as in the proof of
Theorem~\ref{Theorem 3.5}, consider the set
$$T=\bigcup_{1\le i,j\le3}H_{ij}\cup\{0\}$$
where $H_{ij}$ are defined via \eqref{C3 is in M2}. Obviously,
$H_{ij}^*=H_{ji}$ whence $\Tc=\langle T,\cdot,{}^\dag,{}^*\rangle$
is a substructure of $\Ac$ and the partition $\Hc$ of $T$ into the
classes $H_{ij}$ and $\{0\}$ is a congruence on this substructure.
The quotient $\Tc/\Hc$ is then isomorphic to the semigroup
$\mathcal{K}_3$ endowed twice with its unary operation. We
conclude that  $\mathcal{K}_3$ treated this way also belongs to
$\var\Ac$. By Corollary~\ref{identity for rank 1} the identity
$x^2yx=xyx^2$ holds in $\H(\Ac)$ (by which we mean the
substructure of $\Ac$ generated by all elements of the form
$AA^\dag$). Now construct the semigroups $\Tc_k$ (by use of the
identity $x^2yx=xyx^2$) as in Step~2 in the proof of
Theorem~\ref{Theorem 2.1} and endow each of them twice with its
unary operation. The arguments in Steps~3 and~4 in the proof then
show that $\Tc_k$ does not belong to $\var \Ac$ while each
$k$-generated substructure of $\Tc_k$ does belong to $\var\Ac$.
Thus, $\Tc_k$ can play the role of critical structures for
$\var\Ac$ whence the desired conclusion follows by the reasoning
as in Step~1 in the proof of Theorem~\ref{Theorem 2.1}.
\end{proof}

Also in this setting, our result gives rise to a natural question.
\begin{Problem}
\label{problem on Moore-Penrose + conjugation} Is the algebraic
structure
$\langle\mathrm{M}_n(\mathcal{K}),\cdot,{}^\dag,{}^*\rangle$ of
type $(2,1,1)$ not finitely based for each subfield $\mathcal{K}$
of $\bb C$ closed under complex conjugation and for all $n>2$?
\end{Problem}

Here an observation similar to Remark~\ref{rank 1 + invertible}
can be stated: for each con\-jugation-closed subfield
$\mathcal{K}\subseteq\bb C$ and for all $n>2$, the algebraic
structure
$\langle\mathrm{L}_n(\mathcal{K})\cup\mathrm{GL}_n(\mathcal{K}),\cdot,{}^\dag,{}^*\rangle$
consisting of all matrices of rank at most~1 and all invertible
matrices has no finite identity basis.

\subsection{Matrix semigroups with transposition over infinite fields}
\label{usual transposition:infinite fields}

Here we show that the involutory semigroup
$\langle\mathrm{M}_n(\mathcal{K}),\cdot,{}^T\rangle$ is finitely based for any
infinite field $\mathcal{K}$. More precisely, we verify that all identities
holding in $\langle\mathrm{M}_n(\mathcal{K}),\cdot,{}^T\rangle$ follow from the
associativity and the involution laws $(xy)^T = y^Tx^T$, $(x^T)^T= x$. This is
an involutory analogue of a result from \cite{GoMi78} mentioned in the
introduction; the proof given there does not immediately show the intended
analogue, but the ideas below are inspired by the arguments in \cite{GoMi78}.

Let us start with an auxiliary construction and consider first, for an
arbitrary field $\mathcal{K}$, the set $\mathrm{M}_2(\mathcal{K}[x])$ of all
$2\times 2$-matrices over the polynomial ring $\mathcal{K}[x]$. A matrix
$$\begin{pmatrix}
p_{11} & p_{12}\\ p_{21}& p_{22} \end{pmatrix}$$ with $p_{ij}$ non-zero members
of $\mathcal{K}[x]$ is said to be \emph{ascending} if
$$\deg(p_{11})<\min\{\deg(p_{12}),\deg(p_{21})\}, \ \max\{\deg(p_{12}),\deg(p_{21})\}<\deg(p_{22})$$
and \emph{descending} if
$$\deg(p_{11})>\max\{\deg(p_{12}),\deg(p_{21})\}, \ \min\{\deg(p_{12}),\deg(p_{21})\}>\deg(p_{22}.)$$
Let $\mathrm{Asc}$ and $\mathrm{Desc}$ stand for the sets of all ascending and
descending matrices respectively. It is straightforward to see that
$\mathrm{Asc}$ and $\mathrm{Desc}$ are disjoint and both of them are closed
under multiplication and transposition. Set
$$A=\begin{pmatrix} 1 & 0\\ x^2 &x\end{pmatrix};$$ by induction one
obtains that
$$A^n=\begin{pmatrix} 1 & 0 \\p_{n+1}&p_n\end{pmatrix}$$ for polynomials
$p_n$ and $p_{n+1}$ with $\deg(p_{n+1})=n+1$  and $\deg(p_n)=n$. Further, set
$B=A^T$; then for any $n,m\ge 1$ one gets
$$A^nB^m=\begin{pmatrix} 1 & s_{12}\\ s_{21} & s_{22}\end{pmatrix}$$ where
$\deg(s_{12})=m+1$, $\deg(s_{21})=n+1$ and $\deg(s_{22})=m+n+2$, while
$$B^mA^n=\begin{pmatrix} r_{11} & r_{12} \\ r_{21} &
r_{22}\end{pmatrix}$$ where $\deg(r_{11})=m+n+2$,
$\deg(r_{12})=\deg(r_{21})=m+n+1$ and $\deg(r_{22})=m+n$. In particular,
$A^nB^m\in \mathrm{Asc}$ while $B^mA^n\in\mathrm{Desc}$.
\begin{Lemma}\label{lemma:disjoint} For all $m,n\ge 1$, the sets $\mathrm{Asc}\cdot A^n$ and
$\mathrm{Desc}\cdot B^m$ are disjoint.
\end{Lemma}
\begin{proof} Let $a,b,c,d,p,q$ be non-zero polynomials with
$$\deg(a)<\min\{\deg(b),\deg(c)\},\ \max\{\deg(b),\deg(c)\}<\deg(d) \
\text{ and } \ \deg(p)<\deg(q).$$
Then $\left(\begin{smallmatrix} a& b\\ c &
d\end{smallmatrix}\right)$ is a typical matrix in $\mathrm{Asc}$,
$\left(\begin{smallmatrix} d& b\\ c & a\end{smallmatrix}\right)$ is such in
$\mathrm{Desc}$, $A^n$ is of the form $\left(\begin{smallmatrix} 1 & 0\\q  &
p\end{smallmatrix}\right)$ and $B^m$ is of the form $\left(\begin{smallmatrix}
1 & q\\0 & p\end{smallmatrix}\right)$. Now
$$\begin{pmatrix} a & b\\c & d\end{pmatrix} \begin{pmatrix}1 &0\\ q
& p\end{pmatrix} =\begin{pmatrix} a+bq& bp\\ c+dq& dp\end{pmatrix}.$$ Since
$\deg(a)<\deg(b)$ and $\deg(p)<\deg(q)$ we obtain
$$\deg(a+bq)=\deg(bq)> \deg(bp).$$
In particular, for any $C=(c_{ij})\in\mathrm{Asc}\cdot A^n$ we get
\begin{equation}\label{A}
\deg(c_{11})>\deg(c_{12})
\end{equation}

On the other hand,
$$\begin{pmatrix} d & b\\c & a\end{pmatrix} \begin{pmatrix}1 &q\\
0 & p\end{pmatrix} =\begin{pmatrix} d& dq+bp\\
c&cq+ap\end{pmatrix}.$$ Here we have that
$$\deg(d)<\deg(dq)=
\deg(dq+bp).$$  Again, this shows that for any $D=(d_{ij})\in\mathrm{Desc}\cdot
B^m$,
\begin{equation}\label{B}
\deg(d_{11})<\deg(d_{12}).
\end{equation}
Conditions \eqref{A} and \eqref{B} immediately imply that $\mathrm{Asc}\cdot
A^n\cap \mathrm{Desc}\cdot B^m=\varnothing$.
\end{proof}
We are able to draw our first conclusion.
\begin{Prop}
\label{freesubsemigroup} Let $u=u(a,b),v=v(a,b)\in\{a,b\}^+$ be two distinct
words, $A=\left(\begin{smallmatrix} 1 &0\\ x^2 & x\end{smallmatrix}\right)\in
\mathrm{M}_2(\mathcal{K}[x])$ and $B=A^T$. Then $u(A,B)\ne v(A,B)$.
\end{Prop}

\begin{proof}
Suppose that the equality $u(A,B)=v(A,B)$ holds. The matrices $A$ and $B$ are
invertible as  matrices over the field of rational functions over $\mathcal{K}$
whence we can cancel the longest common prefix and the longest common suffix of
the words $u$ and $v$ obtaining new words $u'$ and $v'$ that still fulfil
$u'(A,B)=v'(A,B)$. Observe that none of the words $u'$ and $v'$ are empty
because otherwise $A$ and/or $B$ would be invertible over the polynomial ring
$\mathcal{K}[x]$ which is not true. Thus, we may assume that $u$ and $v$ start
and end with different symbols and this means that the words are of either of
the following two forms:
\begin{gather}\label{eqI}
u=a^{n_1}b^{m_1}\cdots a^{n_s}b^{m_s} \ \text{ and } \
v=b^{\ell_1}a^{k_1}\cdots
b^{\ell_t}a^{k_t}\\
\label{eqII} u=a^{n_1}b^{m_1}\cdots a^{n_{s-1}}b^{m_{s-1}}a^{n_s} \ \text{ and}
\ v=b^{\ell_1}a^{k_1}\cdots b^{\ell_{t-1}}a^{k_{t-1}}b^{k_t}
\end{gather}
with $s,t\ge 1$ and all $n_i,m_i,k_i,\ell_i\ge 1$. In case \eqref{eqI},
$u(A,B)\in \mathrm{Asc}$ and $v(A,B)\in \mathrm{Desc}$ while in case
\eqref{eqII} $u(A,B)\in \mathrm{Asc}\cdot A^{n_s}$ and
$v(A,B)\in\mathrm{Desc}\cdot B^{k_t}$. In any case, using Lemma
\ref{lemma:disjoint} we obtain $u(A,B)\ne v(A,B)$, a contradiction.
\end{proof}

\begin{Thm}
For each infinite field $\mathcal{K}$, all the identities of the involutory
semigroup $\langle\mathrm{M}_n(\mathcal{K}),\cdot,{}^T\rangle$ follow from the
laws $(xy)z=x(yz)$, $(xy)^T = y^Tx^T$, $(x^T)^T= x$.
\end{Thm}

\begin{proof}
If we call identities that follow from the associativity and the involution
laws \emph{trivial}, our claim amounts to saying that
$\langle\mathrm{M}_n(\mathcal{K}),\cdot,{}^T\rangle$ does not satisfy any
non-trivial involutory identity. We show this already holds for the involutory
semigroup $\mathcal{GL}_2(\mathcal{K})=\langle
\mathrm{GL}_2(\mathcal{K}),\cdot,{}^T\rangle$ that obviously embeds into
$\langle\mathrm{M}_n(\mathcal{K}),\cdot,{}^T\rangle$ for each $n$.

It is known (and easy to verify) that the free involutory semigroup on one
generator $z$, say, contains as a unary subsemigroup a free involutory
semigroup on countably many generators, namely, $\FI(Z)$ where
$$Z=\{zz^*z,\ z(z^*)^2z,\ \dots,\ z(z^*)^nz,\ \dots\}.$$
Therefore we only need to verify that $\mathcal{GL}_2(\mathcal{K})$ satisfies
no non-trivial involutory identity in one letter $z$. Such an identity can be
written as $u(z,z^T)=v(z,z^T)$ with $u$ and $v$ being two distinct words.
Proposition~\ref{freesubsemigroup} implies that $u(A,B)\ne v(A,B)$ in
$\mathrm{M}_2(\mathcal{K}[x])$. Setting
$$u(A,B)=\begin{pmatrix} u_{11}&u_{12}\\u_{21} &
u_{22}\end{pmatrix}\ \text{ and } \ v(A,B)=\begin{pmatrix} v_{11} & v_{12}\\
v_{21} & v_{22}\end{pmatrix},$$ we see that for some indices $i$ and $j$ the
polynomials $u_{ij}$ and $v_{ij}$ are distinct whence $u_{ij}-v_{ij}$ is a
non-zero polynomial. Now take
any $\lambda\in \mathcal{K}$ and set $z(\lambda)=\left(\begin{smallmatrix} 1 & 0\\
\lambda^2 & \lambda\end{smallmatrix}\right)$. If the equality
\begin{equation}\label{lambda}
u(z(\lambda),z(\lambda)^T)=v(z(\lambda),z(\lambda)^T)
\end{equation} holds then $\lambda$ must be a root of the non-zero polynomial
$u_{ij}-v_{ij}$ whence (\ref{lambda}) can hold only for finitely many elements
$\lambda$ of $\mathcal{K}$. Since $\mathcal{K}$ is infinite, equality
\eqref{lambda} fails for all but finitely many $\lambda$, and so, in
particular, the identity $u(z,z^T)=v(z,z^T)$ fails in
$\mathcal{GL}_2(\mathcal{K})$.
\end{proof}

Similarly, it can be shown that the involutory semigroup
$\langle\mathrm{M}_n(R),\cdot,{}^*\rangle$ is finitely based for each subring
$R\subseteq \bb C$ closed under complex conjugation---here ${}^*$ stands for
the complex-conjugate transposition $(a_{ij})^*= (\ol{a_{ij}})^T$.  Indeed, we
have already mentioned that the two matrices $\zeta$ and $\eta$ in~\eqref{free
subgroup} generate a free subgroup of $\langle\mathrm{SL}_2(\bb
Z),\cdot,{}^{-1}\rangle$ and hence a free subsemigroup of
$\langle\mathrm{SL}_2(\bb Z), \cdot\rangle$. But $\eta=\zeta^*$ whence the
involutory subsemigroup in $\langle\mathrm{SL}_2(\bb Z),\cdot,{}^*\rangle$
generated by $\zeta$ is isomorphic to the free involutory semigroup
$\FI(\{\zeta\})$.

We note that Theorem~\ref{Theorem 2.1} and Corollary~\ref{identity for rank 1}
prove the non-existence of a finite identity bases for the unary subsemigroup
of $\langle\mathrm{M}_n(\bb C),\cdot,^*\rangle$ [respectively
$\langle\mathrm{M}_n(\bb R),\cdot,{}^T\rangle$] that consists of all matrices
of rank at most~1 together with all unitary [respectively all orthogonal]
matrices.

\subsection{Matrix semigroups with transposition over finite fields}
\label{usual transposition:finite fields}

Now we demonstrate that the case of finite fields is completely opposite to the
one of infinite fields with respect to the finite basis problem for matrix
semigroups with transposition. We start with considering $2\times 2$-matrices.
Here Theorem~\ref{Theorem 2.1} solves the finite basis problem in the negative
for the involutory semigroup
$\langle\mathrm{M}_2(\mathcal{K}),\cdot,{}^T\rangle$ for each finite field
$\mathcal{K}$ except $\mathcal{K}=\mathbb{F}_2$, the $2$-element field. (The
exception is due to the fact that the involutory semigroup
$\langle\mathrm{M}_2(\bb F_2),\cdot,{}^T\rangle$ satisfies the identity
$$(xx^T)^3(yy^T)^3=(yy^T)^3(xx^T)^3$$
which does not hold in $\mathcal{K}_3$; consequently, $\mathcal{K}_3$ is not in
$\var\langle\mathrm{M}_2(\bb F_2),\cdot,{}^T\rangle$ and Theorem~\ref{Theorem
2.1} does not apply here.) In the following theorem, we shall demonstrate the
application of Theorem \ref{Theorem 2.1} only in the case when $\mathcal{K}$
has odd characteristic. With some additional effort we could include also the
case when the characteristic of $\mathcal{K}$ is $2$ and $\vert K\vert\ge 4$.
We shall omit this since that case will be covered by a different kind of proof
later.

\begin{Thm}\label{Theorem 3.6}
For each finite field $\mathcal{K}$ of odd characteristic, the involutory semigroup
$\langle\mathrm{M}_2(\mathcal{K}),\cdot,{}^T\rangle$ has no finite identity basis.
\end{Thm}

\begin{proof}
Let $\Sc=\langle\mathrm{M}_2(\mathcal{K}),\cdot,{}^T\rangle$. As in the proof of Theorem \ref{Theorem 3.5} one shows
that $\mathcal{K}_3$ is in $\var\Sc$.  Furthermore, let $d$ be the exponent of the group $\mathrm{GL}_2(\mathcal{K})$.
By Corollary~\ref{identity for rank 1}, each group in $\P_d(\var\Sc)=\var\P_d(\Sc)$ satisfies the identity
$x^2yx=xyx^2$ and therefore is abelian. On the other hand, as in the proof of Theorem~\ref{both operations}, the group
$\mathcal{G}=\{A\in \mathrm{GL}_2(\mathcal{K})\mid A^T=A^{-1}\}$ is in $\var\Sc$ but is non-abelian. Thus,
Theorem~\ref{Theorem 2.1} applies.
\end{proof}

The next theorem contains the even characteristic case and proves, in fact, a stronger assertion.

\begin{Thm}\label{inherently degree 2}
For each finite field $\mathcal{K}=\langle K,+,\cdot\rangle$ with $\vert K\vert\mathrel{\not\equiv} 3\ (\bmod\ 4)$, the
involutory semigroup $\langle\mathrm{M}_2(\mathcal{K}),\cdot,{}^T\rangle$ is inherently \nfb.
\end{Thm}

\begin{proof}
As mentioned in Remark \ref{remark on MP inverse}, there exists $x\in K$ for which $1+x^2=0$. Now consider the
following matrices:
\begin{gather*}
H_{11}=\begin{pmatrix} 1 & x\\ x &x^2\end{pmatrix},\ H_{12}=\begin{pmatrix} 1 & 0\\ x &0\end{pmatrix},\
H_{21}=\begin{pmatrix} 1 & x\\ 0 & 0\end{pmatrix},\
H_{22}=\begin{pmatrix} 1 & 0\\ 0 & 0\end{pmatrix},\\
I=\begin{pmatrix}1 & 0\\ 0& 1\end{pmatrix},\ O=\begin{pmatrix} 0 & 0 \\ 0 & 0 \end{pmatrix}.
\end{gather*}
Then the set $M=\{H_{11},H_{12},H_{21}, H_{22},I, O\}$ is closed under multiplication and transposition, hence
$\mathcal{M}=\langle M, \cdot,{}^T\rangle$ is an involutory subsemigroup of
$\langle\mathrm{M}_2(\mathcal{K}),\cdot,{}^T\rangle$. The mapping $\TA\to \mathcal{M}$ given by
$$(i,j)\mapsto H_{ij},\ 0\mapsto O,\ 1\mapsto I$$
is an isomorphism of involutory semigroups. The result now follows from Corollary~\ref{twisted A}.
\end{proof}

The case of matrix semigroups of size greater than 2 is similar.

\begin{Thm} \label{inherently degree 3}
For  each finite field $\mathcal{K}$, the involutory semigroup
$\langle\mathrm{M}_n(\mathcal{K}),\cdot,{}^T\rangle$ with $n\ge 3$ is
inherently \nfb.
\end{Thm}

\begin{proof}
It follows from the Chevalley-Warning theorem \cite[Corollary~2 in \S1.2]{Serre} that there exist $x,y\in K$ satisfying
$1+x^2+y^2=0$. Now consider the following matrices:
\begin{gather*}
H_{11}=\begin{pmatrix} 1 & x & y\\ x &x^2& xy\\ y & xy & y^2\end{pmatrix},\ H_{12}=\begin{pmatrix} 1 & 0 & 0\\ x &0
&0\\ y & 0 & 0\end{pmatrix},\
H_{21}=\begin{pmatrix} 1 & x & y\\ 0 & 0 & 0\\ 0 & 0&0\end{pmatrix},\\
H_{22}=\begin{pmatrix} 1 & 0 & 0\\ 0 & 0&0\\0&0&0\end{pmatrix},\ I=\begin{pmatrix} 1&0&0\\0&1&0\\0&0&1\end{pmatrix},\
O=\begin{pmatrix} 0&0&0\\0&0&0\\0&0&0\end{pmatrix}.
\end{gather*}
Again, the set $M=\{H_{11},H_{12},H_{21},H_{22},I,O\}$ is closed under multiplication and transposition, and as in the
previous proof, $\mathcal{M}=\langle M ,\cdot,{}^T\rangle$ forms an involutory subsemigroup of
$\langle\mathrm{M}_3(\mathcal{K}),\cdot,{}^T\rangle$ that is isomorphic with $\TA$. Hence $\langle
\mathrm{M}_3(\mathcal{K}),\cdot,{}^T\rangle$ is inherently \nfb. The assertion for $\mathrm{M}_n(\mathcal{K})$ for
$n\ge 3$ now follows in an obvious way.
\end{proof}

\begin{Rmk}
The statements of Theorems \ref{Theorem 3.6}, \ref{inherently degree 2}, \ref{inherently degree 3} remain valid if the
unary operation $A\mapsto A^T$ is replaced with an operation of the form $A\mapsto A^{\si T}$ for any automorphism
$\si$ of $\mathcal{K}$, where $(a_{ij})^{\si T}:=(a_{ij}^\si)^T$.
\end{Rmk}

We are ready to prove the main result of this subsection. With the exception of the `only if' part of item (2), this is
a summary of Theorems \ref{Theorem 3.6}, \ref{inherently degree 2}, and \ref{inherently degree 3}.
\begin{Thm}\label{main result matrix involution} Let
$\mathcal{K}=\langle K,+,\cdot\rangle$ be a finite field. Then
\begin{enumerate}
\item the involutory semigroup $\langle \mathrm{M}_n(\mathcal{K}),\cdot,{}^T\rangle$
is not finitely based;
\item the involutory semigroup $\langle \mathrm{M}_n(\mathcal{K}),\cdot,{}^T\rangle$
is inherently nonfinitely ba\-sed if and only if either $n\ge 3$ or $n=2$ and $\vert K\vert\mathrel{\not\equiv 3}
(\bmod\ 4)$.
\end{enumerate}
\end{Thm}

\begin{proof}
The only assertion of this theorem not covered by our previous results is that
$\langle\mathrm{M}_2(\mathcal{K}),\cdot,{}^T\rangle$ is \textbf{not} inherently
\nfb\ if $\vert K\vert\equiv 3\ (\bmod\ 4)$. We employ Proposition~\ref{NINFB}
to prove this.

Recall that the condition $\vert K\vert\equiv 3\ (\bmod\ 4)$ corresponds precisely to the case when each matrix $A$ in
$\langle\mathrm{M}_2(\mathcal{K}),\cdot,{}^T\rangle$ admits a Moore-Penrose inverse $A^\dag$ (Remark~\ref{remark on MP
inverse}). Let $A$ be a matrix of rank 1; by (\ref{MP is scalar multiple}) there exists a scalar $\alpha\in
\mathcal{K}\setminus\{0\}$ such that $\alpha A^\dag= A^T$. Let $r=\vert K\vert -1$; then $\alpha^r=1$. Since the
multiplicative subgroup of $\mathcal{K}$ is a cyclic subgroup of $\mathrm{GL}_2(\mathcal{K})$, the number $r$ divides
the exponent $d$ of $\mathrm{GL}_2(\mathcal{K})$ whence $\alpha^d=1$. Consequently,
$$A(A^TA)^d= A(\alpha A^\dag A)^d=\alpha^d A(A^\dag A)^d=A.$$
If $A\in\mathrm{GL}_2(\mathcal{K})$, we also have $A=A(A^TA)^d$ because $(A^TA)^d$ is the identity matrix; clearly, the
equality $A=A(A^TA)^d$ holds also for the case when $A$ is the zero matrix. Summarizing, we conclude that the identity
$x=x(x^Tx)^d$ holds in the involutory \sm\ $\langle\mathrm{M}_2(\mathcal{K}),\cdot,{}^T\rangle$. Setting
$\om(x):=x^T(xx^T)^{d-1}$, we see that $\langle\mathrm{M}_2(\mathcal{K}),\cdot,{}^T\rangle$ satisfies the identity
$x=x\om(x)x$, as required by Proposition~\ref{NINFB}.
\end{proof}

\begin{Rmk}
It is known~\cite[Corollary~6.2]{sapirburnside} that the matrix semigroup
$\langle\mathrm{M}_n(\mathcal{K}),\cdot\rangle$ is inherently \nfb\ (as a plain \sm) for every finite field
$\mathcal{K}$. Thus, the involutory semigroups $\langle\mathrm{M}_2(\mathcal{K}),\cdot,{}^T\rangle$ over finite fields
$\mathcal{K}$ such that $\vert K\vert\equiv 3\ (\bmod\ 4)$ provide a natural series of unary semigroups whose
equational properties essentially differ from the equational properties of their semigroup reducts.
\end{Rmk}

\subsection{Matrix semigroups with symplectic transpose}
\label{symplectic case}

For a $2m\times 2m$-matrix
$$X=\begin{pmatrix} A & B\\ C & D\end{pmatrix}$$
with $A,B,C,D$ being $m\times m$-matrices over any field $\mathcal{K}$, the \emph{symplectic transpose} $X^S$ is
defined by
$$X^S=\begin{pmatrix} D^T & -B^T\\ -C^T & A^T\end{pmatrix},$$
see, e.g., \cite[(5.1.1)]{Procesi}. The symplectic transpose is an involution of
$\left<\mathrm{M}_{2m}(\mathcal{K}),\cdot\right>$ whose properties essentially differ from those of the usual
transposition. In fact, every involution of the semigroup $\left<\mathrm{M}_n(\mathcal{K}),\cdot\right>$ that fixes the
scalar matrices is in a certain sense similar to either the usual transposition or the symplectic transpose\footnote{We
do not want to formalize this phrase in general because its meaning actually depends on some intrinsic properties of
the field $\mathcal{K}$. In the simplest case, when $\mathcal{K}$ is algebraically closed of characteristic $\ne 2$,
the similarity takes the strongest form: every involutory semigroup of the form
$\langle\mathrm{M}_n(\mathcal{K}),\cdot,{}^*\rangle$ such that ${}^*$ fixes the scalar matrices is isomorphic to either
$\langle\mathrm{M}_n(\mathcal{K}),\cdot,{}^T\rangle$ or $\langle\mathrm{M}_n(\mathcal{K}),\cdot,{}^S\rangle$ (in the
latter case $n$ should be even). This is well known for involutions that respect the addition of matrices (see, e.g.,
\cite[Corollary 14.2]{Procesi-INV}) but it easily follows from a classic result by Khalezov~\cite{Kha54} that every
involution of the semigroup $\left<\mathrm{M}_n(\mathcal{K}),\cdot\right>$ automatically preserves addition.}.

The definition of the symplectic transpose resembles that of the
involution in the twisted Brandt monoid $\mathcal{TB}^1_2$
(defined in terms of $2\times 2$-matrices) and leads to the
following application.
\begin{Thm}\label{symplecticinvolution}
The involutory semigroup $\left<\mathrm{M}_{2m}(\mathcal{K}),\cdot,{}^S\right>$ is inherently \nfb\ for each $m\ge 1$
and each finite field $\mathcal{K}=\left<K,+,\cdot\right>$.
\end{Thm}
\begin{proof} Consider the following sets of $2m\times 2m$-matrices:
\begin{gather*}
H_{11}=\left\{\pm \begin{pmatrix}O_m & I_m \\ O_m&O_m
\end{pmatrix}\right\},\
H_{12}=\left\{\pm \begin{pmatrix}I_m & O_m \\ O_m&O_m \end{pmatrix}\right\},\\
H_{21}=\left\{\pm \begin{pmatrix}O_m & O_m \\ O_m&I_m
\end{pmatrix}\right\},\
H_{22}=\left\{\pm \begin{pmatrix}O_m & O_m \\ I_m&O_m
\end{pmatrix}\right\}
\end{gather*}
where for any positive integer $k$, we denote be $I_k$, respectively, $O_k$ the identity, respectively, zero $k\times
k$-matrix. Let
$$T=\bigcup_{1\le i,j\le 2}H_{ij}\cup \{O_{2m}, I_{2m}\}.$$
The set $T$ is closed under multiplication and symplectic transposition whence
$\left<T,\cdot,{}^S\right>$ forms an involutory subsemigroup of
$\left<\mathrm{M}_{2m}(\mathcal{K}),\cdot,{}^S\right>$. On the other hand, the
mapping
$$H_{ij}\mapsto (i,j),\ I_{2m}\mapsto 1,\ O_{2m} \mapsto 0$$
is a homomorphism of $\left<T,\cdot,{}^S\right>$ onto $\mathcal{TB}^1_2$.
Altogether, the twisted Brandt monoid $\mathcal{TB}^1_2$ is a homomorphic image
of an involutory subsemigroup of
$\left<\mathrm{M}_{2m}(\mathcal{K}),\cdot,{}^S\right>$.
\end{proof}

\begin{Rmk}
It is easy to see that the involutory semigroup
$\langle\mathrm{M}_m(\mathcal{K}),\cdot,{}^T\rangle$ embeds into
$\left<\mathrm{M}_{2m}(\mathcal{K}),\cdot,{}^S\right>$ via the mapping
$A\mapsto\begin{pmatrix} A & O_m\\ O_m & A\end{pmatrix}$. Therefore, if $m>1$
and $\mathcal{K}$ is an infinite field, then
$\left<\mathrm{M}_{2m}(\mathcal{K}),\cdot,{}^S\right>$ satisfies no non-trivial
involutory identity and is finitely based (see Subsection~\ref{usual
transposition:infinite fields}). The $2\times 2$-matrices over an infinite
field fulfill non-trivial identities involving multiplication and the
symplectic transpose, for instance, $xx^Sy=yxx^S$ or $xx^S=x^Sx$. In fact, we
have verified that these two identities together with the associativity and the
involution laws $(xy)^S = y^Sx^S$, $(x^S)^S= x$ form an identity basis for
$\left<\mathrm{M}_{2}(\mathcal{K}),\cdot,{}^S\right>$ with $\mathcal{K}$
infinite. (The proof of this will be published separately.) Thus, the
involutory semigroup $\left<\mathrm{M}_{2m}(\mathcal{K}),\cdot,{}^S\right>$ is
finitely based if and only if $\mathcal{K}$ is an infinite field.
\end{Rmk}

\subsection{Boolean matrices}
\label{Boolean case}

Recall that a \emph{Boolean} matrix is a matrix with entries $0$
and $1$ only. The multiplication of such matrices is as usual,
except that addition and multiplication of the entries is defined
as: $a+b=\max\{a,b\}$ and $a\cdot b=\min\{a,b\}$. Let $B_n$ denote
the set of all Boolean $n\times n$-matrices. It is well known that
the \sm\ $\langle B_n,\cdot\rangle$ is essentially the same as the
semigroup of all binary relations on an $n$-element set subject to
the usual composition of binary relations. The operation ${}^T$ of
forming the matrix transpose then corresponds to the operation of
forming the dual binary relation.

\begin{Thm}
\label{Theorem 3.7}
The involutory semigroup $\Bc_n=\langle B_n,\cdot,{}^T\rangle$
of all Boolean $n\times n$-matrices endowed with transposition
is inherently \nfb.
\end{Thm}

\begin{proof} Consider the Boolean matrices
\begin{gather*}
B_{11}=\begin{pmatrix} 0&1\\1&1\end{pmatrix},\
B_{12}=\begin{pmatrix} 1&0\\1&1 \end{pmatrix},\
B_{21}=\begin{pmatrix} 1&1\\0&1 \end{pmatrix},\
B_{22}=\begin{pmatrix} 1&1\\1&0\end{pmatrix},\\
O=\begin{pmatrix} 1&1\\1&1 \end{pmatrix},\
I=\begin{pmatrix}1&0\\0&1 \end{pmatrix}.
\end{gather*}
The set $M=\{B_{11}, B_{12}, B_{21}, B_{22},O,I\}$ is closed under
multiplication and transposition whence $\mathcal{M}=\langle
M,\cdot,{}^T\rangle$ is an involutory subsemigroup of $\Bc_2$. The
mapping $\TB\to \mathcal{M}$ given by
$$ (i,j)\mapsto B_{ij},\ 0\mapsto O,\ 1\mapsto I$$
is an isomorphism of involutory semigroups. By
Corollary~\ref{twisted Brandt} $\mathcal{M}$ is inherently
\nfb\ whence so is $\Bc_2$. Since $\Bc_2$ can be embedded
as an involutory semigroup into $\Bc_n$ for each $n$,
the result follows.
\end{proof}

\begin{Rmk}
Our proof of Theorem~\ref{Theorem 3.7} also applies to some
important involutory subsemigroups of $\Bc_n$. In order to
introduce an interesting instance, recall that Boolean
$n\times n$-matrices are in a 1-1 correspondence with
bipartite directed graphs whose parts are of size $n$:
the bipartite graph of a matrix $A=(a_{ij})$ has the row
set and the column set of $A$ as its parts and has an edge
from the $i^\mathrm{th}$ row to the $j^\mathrm{th}$ column
if and only if $a_{ij}=1$. If the graph of $A$ admits
a perfect matching (i.e.\ a set of edges so that every
vertex is incident to precisely one of them), $A$ is said
to be a \emph{Hall matrix} (the name suggested in~\cite{Kim}
is, of course, inspired by Hall's marriage theorem). It
is easy to see that the collection $H\!B_n$ of all Hall
$n\times n$-matrices is closed under multiplication and
transposition. Since all the matrices $B_{11}, B_{12}, B_{21},
B_{22},O,I$ from the above proof are Hall matrices, we readily
conclude that the involutory semigroup $\langle H\!B_n,\cdot,{}^T\rangle$
is inherently \nfb.
\end{Rmk}

\begin{Rmk}
We can unify Theorem~\ref{Theorem 3.7} and some results in
Subsection~\ref{usual transposition:finite fields} by considering matrices
over \textbf{semirings}. A \emph{semiring} is an algebraic
structure $\mathcal{L}=\langle L,+,\cdot\rangle$ of type $(2,2)$
such that $\langle L,+\rangle$ is a commutative semigroup,
$\langle L,\cdot\rangle$ is a semigroup and multiplication
distributes over addition. From the proofs of
Theorems~\ref{inherently degree 2} and~\ref{Theorem 3.7}
we see that the involutory matrix semigroup
$\langle\mathrm{M}_n(\mathcal{L}),\cdot,{}^T\rangle$ over a finite
semiring is inherently \nfb\ whenever
the semiring $\mathcal{L}$ has a zero 0 (that is, a neutral
element for $\langle L,+\rangle$ which is at the same time an
absorbing element for $\langle L,\cdot\rangle$) and satisfies
one of the following two conditions:
\begin{enumerate}
\item there exist (not necessarily distinct) elements $e,x\ne 0$
such that $e^2=e$, $ex=xe=x$, $e+x^2=0$;
\item there exists an element $e\ne 0$ such that $e^2=e=e+e$.
\end{enumerate}
We have already met an infinite series of semirings satisfying
(1): it consists of the finite fields $\mathcal{K}=\langle
K,+,\cdot\rangle$ with $\vert K\vert\mathrel{\not\equiv 3}(\bmod\
4)$. It should be noted that semirings satisfying (2) are even
more plentiful: for example, finite distributive lattices as well
as the power semirings of finite semigroups (with the subset union
as addition and the subset product as multiplication) fall in this
class.
\end{Rmk}

A Boolean matrix $A=(a_{ij})$ is said to be \emph{upper triangular}
if $a_{ij}=0$ whenever $i>j$. Let $BT_n$ stand for the set of all
Boolean upper triangular $n\times n$-matrices. The semigroups
$\langle BT_n,\cdot\rangle$ play an important role in the
theory of formal languages, see~\cite{PiSt85}; their identities
have been studied by the third author and Goldberg in~\cite{VoGo04}.
Observe that this \sm\ admits quite a natural unary operation:
the reflection with respect to the secondary diagonal (the diagonal
from the top right to the bottom left corner). We denote by $A^D$
the result of applying this operation to the matrix $A$. It is easy
to see that the operation $A\mapsto A^D$ is in fact an involution;
this follows, for instance, from the fact that $A^D=JA^TJ$ where $J$
is the Boolean matrix with 1s in the secondary diagonal and 0s elsewhere.

\begin{Thm}
\label{Theorem 3.8} For each integer $n\ge 3$, the involutory
semigroup $\BT_n=\langle BT_n,\cdot,{}^D\rangle$ of all Boolean
upper triangular $n\times n$-matrices endowed with the reflection
with respect to the secondary diagonal is inherently \nfb.
\end{Thm}

\begin{proof}
Consider the involutory submonoid $\mathcal{M}$ in $\BT_n$
generated by the following two Boolean matrices:
$$X=\begin{pmatrix}
1 & 0 & \dots & 0 & 0 \\
0 & 0 & \dots & 0 & 0\\
\vdots & \vdots & \ddots & \vdots & \vdots\\
0 & 0 & \dots & 0 & 0\\
0 & 0 & \dots & 0 & 1
\end{pmatrix}\quad
\text{and}\quad
Y=\begin{pmatrix}
1 & 1 & \dots & 1 & 0 \\
0 & 0 & \dots & 0 & 1\\
\vdots & \vdots & \ddots & \vdots & \vdots\\
0 & 0 & \dots & 0 & 1\\
0 & 0 & \dots & 0 & 1
\end{pmatrix}. $$
Clearly, for each matrix $(m_{ij})$ in this submonoid one has
$m_{11}=m_{nn}=1$, whence the set of all matrices $(m_{ij})$ such
that $m_{1n}=1$ forms an ideal in $\mathcal{M}$. We denote this
ideal by $\mathcal{N}$. A straightforward calculation shows that,
besides $X$, $Y$, and the identity matrix $I$, only the two
matrices
$$XY=\begin{pmatrix}
1 & 1 & \dots & 1 & 0 \\
0 & 0 & \dots & 0 & 0\\
\vdots & \vdots & \ddots & \vdots & \vdots\\
0 & 0 & \dots & 0 & 0\\
0 & 0 & \dots & 0 & 1
\end{pmatrix}\quad \text{and}\quad
YX=\begin{pmatrix}
1 & 0 & \dots & 0 & 0 \\
0 & 0 & \dots & 0 & 1\\
\vdots & \vdots & \ddots & \vdots & \vdots\\
0 & 0 & \dots & 0 & 1\\
0 & 0 & \dots & 0 & 1
\end{pmatrix}$$ belong to
$\mathcal{M}\setminus\mathcal{N}$. This allows one to organize the
following bijection between $\mathcal{M}\setminus\mathcal{N}$ and
the set of non-zero matrices in $\TA$:
$$I\mapsto\left(\begin{smallmatrix}
1 & 0 \\ 0 & 1
\end{smallmatrix}\right),\quad
X\mapsto\left(\begin{smallmatrix}
1 & 0 \\ 1 & 0
\end{smallmatrix}\right),\quad
Y\mapsto\left(\begin{smallmatrix}
0 & 1 \\ 0 & 0
\end{smallmatrix}\right),\quad
XY\mapsto\left(\begin{smallmatrix}
0 & 1 \\ 0 & 1
\end{smallmatrix}\right),\quad
YX\mapsto\left(\begin{smallmatrix}
1 & 0 \\ 0 & 0
\end{smallmatrix}\right). $$
One easily checks that extending this bijection to $\mathcal{M}$
by sending all elements from $\mathcal{N}$ to
$\left(\begin{smallmatrix}
0 & 0 \\ 0 & 0
\end{smallmatrix}\right)$ yields an involutory \sm\
homomorphism from $\mathcal{M}$ onto $\TA$. Thus, $\TA$ as
a homomorphic image of an involutory sub\sm\ in $\BT_n$ belongs to the
involutory \sm\ variety generated by $\BT_n$.
Corollary~\ref{twisted A} implies that $\BT_n$ is inherently
\nfb.
\end{proof}

\begin{Rmk}
In \cite{VoGo04} it was shown that the \sm\ $\langle BT_n,\cdot\rangle$
is inherently \nfb\ for $n\ge4$. However, the construction
used there does not imply the same fact for the involutory case. Our
proof of Theorem~\ref{Theorem 3.8} follows a different construction
suggested for the plain \sm\ case by Li and Luo~\cite{LL}. Li and Luo
have also verified that the \sm\ $\langle BT_2,\cdot\rangle$ is finitely
based. We do not know whether or not the involutory semigroup $\BT_2$
is finitely based.
\end{Rmk}

Other interesting involutory semigroups of Boolean matrices include
the \sm\ $\BR_n=\langle B\!R_n,\cdot,{}^T\rangle$ of all Boolean
$n\times n$-matrices with 1s on the main diagonal (such matrices
correspond to reflexive binary relations) and the \sm\ $\BU_n=\langle
BU_n,\cdot,{}^D\rangle$ of all Boolean upper triangular
$n\times n$-matrices with 1s on the main diagonal. (The unary operation
is the usual transpose in the former case and the reflection
with respect to the secondary diagonal in the latter one.) Our present
methods do not yet suffice to handle the finite basis problem for
these unary semigroups so we just formulate
\begin{Problem}
Are the involutory semigroups $\BR_n$ and $\BU_n$ finitely based for $n\ge 3$?
\end{Problem}
The involutory semigroups $\BR_2$ and $\BU_2$ are easily seen to be
finitely based. The finite basis problem for the plain semigroups
$\langle B\!R_n,\cdot\rangle$ and $\langle BU_n,\cdot\rangle$ has been
solved by the third author~\cite{Vo04}.

\medskip

\noindent\textbf{Acknowledgement.} The second author was supported by Grant No.\ 174019 of the Ministry of Science and
Technological Development of the Republic of Serbia. The third author acknowledges support from the Ministry for Education
and Science of Russia, grant 2.1.1/13995, and from the Russian Foundation for Basic Research, grant 10-01-00524.


\begin{thebibliography}{99}
\frenchspacing

\bibitem{Ami74}
Amitsur, S. A.: Polynomial identities. Israel J. Math. \textbf{19}, 183--199 (1974)

\bibitem{AM}
Ara\'ujo, J., Mitchell, J. D.: An elementary proof that every singular matrix is a product of idempotent matrices.
Amer.\ Math.\ Monthly \textbf{112}, 641--645 (2005)

\bibitem{A1}
Auinger, K.: Strict regular $*$-semigroups. In: Proceedings of the Conference on Semigroups with Applications, Howie,
J. M., Munn, W. D., Weinert, H.-J. (eds.), pp.\ 190--204,  World Scientific, Singapore (1992)

\bibitem{BIG}
Ben-Israel, A., Greville, Th.: Generalized Inverses: Theory and Applications. Sprin\-ger-Verlag,
Berlin--Heidelberg--New York (2003)

\bibitem{BuSa81}
Burris, S., Sankappanavar, H. P.: A Course in Universal Algebra. Springer-Verlag, Berlin--Heidelberg--New York (1981)

\bibitem{CP}
Clifford, A. H., Preston, G. B.: The Algebraic Theory of Semigroups, Vol. I. Amer.\ Math.\ Soc., Providence (1961)

\bibitem{Cline}
Cline, R. E.: Note on the generalized inverse of the product of matrices. SIAM Review \textbf{6}, 57--58 (1964)

\bibitem{CoKo04}
Colombo, J., Koshlukov, P.: Central polynomials in the matrix algebra of order two. Linear Algebra Appl. \textbf{377},
53--67 (2004)

\bibitem{DaRa99}
D'Amour, A., Racine, M.: *-Polynomial identites of matrices with the transpose involution: the low degrees. Trans.\
Amer.\ Math.\ Soc. \textbf{351}, 5089--5106 (1999)

\bibitem{DaRa04}
D'Amour, A., Racine, M.: *-Polynomial identites of matrices with the symplectic involution: the low degrees. Comm.\
Algebra \textbf{32}, 895–-918 (2004)

\bibitem{Dolinka}
Dolinka, I.: On identities of finite involution semigroups. Semigroup Forum \textbf{80}, 105--120 (2010)

\bibitem{D}
Drazin, M. P.: Regular semigroups with involution. In: Proceedings of the Symposium on Regular Semigroups,
pp. 29--46,  Northern Illinois University, De Kalb (1979)

\bibitem{Dr81}
Drensky, V. S.: A minimal basis of identities for a second-order matrix algebra over a field of characteristic $0$.
Algebra i Logika \textbf{20}, 282--290 (1981) [Russian;  English transl. Algebra and Logic \textbf{20}, 188--194
(1982)]

\bibitem{DrFo04}
Drensky, V. S., Formanek, E.: Polynomial Identity Rings. Birkh\"auser, Basel (2004)

\bibitem{DrGi95}
Drensky, V. S., Giambruno, A.: On the *-polynomial identities of minimal degree for matrices with involution. Boll.\
Unione Math.\ Ital. Ser. A \textbf{9}, 471--482 (1995)

\bibitem{Erdos}
Erdos, J. A.: On products of idempotent matrices. Glasgow Math.\ J. \textbf{8}, 118--122 (1967)

\bibitem{Ge81}
Genov, G. K.: A basis of identities of the algebra of third-order matrices over a finite field. Algebra i Logika
\textbf{20}, 385--388 (1981) [Russian;  English transl. Algebra and Logic \textbf{20}, 241--257 (1982)]

\bibitem{GeSi82}
Genov, G. K., Siderov, P. N.: A basis of the identities of the fourth order matrix algebra over a finite field. I, II.
Serdica \textbf{8}, 313--323, 351--366 (1982) [Russian]

\bibitem{GePe85}
Gerhard, J. A., Petrich, M.: Free involutorial completely simple semigroups. Canad. J. Math. \textbf{37}, 271--295
(1985)

\bibitem{Gi90}
Giambruno, A.: On *-polynomial identities for $n\times n$-matrices. J. Algebra \textbf{133}, 433-–438 (1990)

\bibitem{GiZa05}
Giambruno, A., Zaicev, M.: Polynomial Identities and Asymptotic Methods. Amer.\ Math.\ Soc., Providence (2005)

\bibitem{GoMi78}
Golubchik, I. Z., Mikhalev, A. V.: A note on varieties of semiprime rings with semigroup identities. J. Algebra
\textbf{54}, 42--45 (1978)

\bibitem{KaRo05}
Kanel-Belov, A., Rowen, L. H.: Computational Aspects of Polynomial Identities. A~K~Peters Ltd., Wellesley (2005)

\bibitem{Ke87}
Kemer, A. R.: The finite basis property of identities of associative algebras. Algebra i Logika \textbf{26}, 597--641
(1987) [Russian; English transl. Algebra and Logic \textbf{26}, 362--397 (1987)]

\bibitem{Ke91}
Kemer, A. R.: Ideals of Identities of Associative Algebras. Amer.\ Math.\ Soc., Providence (1991)

\bibitem{Kha54}
Khalezov, E. A.: Isomorphisms of matrix semigroups. Ivanov.\ Gos.\ Ped.\ Inst., Uchenye Zap., Fiz.-Mat.\ Nauki
\textbf{5}, 42--56 (1954) [Russian]

\bibitem{Kim}
Kim, K. H.: The semigroups of Hall relations. Semigroup Forum \textbf{9}, 253--260 (1974)

\bibitem{KR}
Kim, K. H., Roush, F.: On groups in varieties of semigroups. Semigroup Forum \textbf{16}, 201--202 (1978)

\bibitem{kleiman}
Kleiman, E. I.: Bases of identities of varieties of inverse semigroups. Sib.\ Mat.\ Zh. \textbf{20}, 760--777 (1979)
[Russian; English transl.\ Sib.\ Math. J. \textbf{20}, 530--543 (1979)]

\bibitem{Ko01}
Koshlukov, P.: Basis of the identities of the matrix algebra of order two over a field of characteristic $p\ne2$. J.
Algebra \textbf{241}, 410--434 (2001)

\bibitem{Kr73}
Kruse, R. L.: Identities satisfied by a finite ring. J. Algebra \textbf{26}, 298--318 (1973)

\bibitem{LL}
Li, J. R., Luo, Y. F.: On the finite basis problem for the monoids of triangular boolean
matrices. Algebra Univers. (to appear)

\bibitem{LidlNiederreiter}
Lidl, R., Niederreiter, H.: Finite Fields. Addison-Wesley, Cambridge (1997)

\bibitem{Lv73}
L'vov, I. V.: Varieties of associative rings. I. Algebra i Logika \textbf{12}, 269--297 (1973) [Russian; English
transl. Algebra and Logic \textbf{12}, 150--167 (1973)]

\bibitem{mks}
Magnus, W., Karras, A., Solitar, D.: Combinatorial Group Theory. Wiley, New York--London--Singapore (1966)

\bibitem{MaKu78}
Mal'tsev, Yu. N., Kuz'min, E. N.: A basis for the identities of the algebra of second-order matrices over a finite
field. Algebra i Logika \textbf{17}, 28--32 (1978) [Russian; English transl.\ Algebra and Logic \textbf{17}, 18--21
(1978)]

\bibitem{margolissapir}
Margolis, S. W., Sapir, M. V.: Quasi-identities of finite semigroups and symbolic dynamics. Israel J. Math.
\textbf{92}, 317--331 (1995)

\bibitem{Meyer}
Meyer, C. D.: Matrix Analysis and Applied Linear Algebra. SIAM, Philadelphia (2000)

\bibitem{moore}
Moore, E. H.: On the reciprocal of the general algebraic matrix. Bull.\ Amer.\ Math.\ Soc. \textbf{26}, 394--395 (1920)

\bibitem{Ne}
Neumann, H.: Varieties of Groups. Springer-Verlag, Berlin--Heidelberg--New York (1967)

\bibitem{oatespowell}
Oates, S., Powell, M. B.: Identical relations in finite groups. J. Algebra \textbf{1}, 11--39 (1964)

\bibitem{P}
Penrose, R.: A generalized inverse for matrices. Proc.\ Cambridge Phil.\ Soc. \textbf{51}, 406--413 (1955)

\bibitem{PiSt85}
Pin, J.-E., Straubing, H.: Monoids of upper triangular matrices. In: Semigroups. Structure and Universal Algebraic
Problems, Poll\'ak, Gy., Schwarz, \v{S}t., Steinfeld, O. (eds.), Colloq. Math. Soc. J\'anos Bolyai \textbf{39}, pp.
259--272, North-Holland, Amsterdam--Oxford--New York (1985)

\bibitem{Procesi-INV}
Procesi, C.: The invariant theory of $n\times n$ matrices. Adv.\ Math. \textbf{19}, 306--381 (1976)

\bibitem{Procesi}
Procesi, C.: Lie Groups: an Approach through Invariants and Representations. Springer-Verlag, Berlin--Heidelberg--New
York (2006)

\bibitem{Ra73}
Razmyslov, Yu. P.: Finite basing of the identities of a matrix algebra of second order over a field of characteristic
zero. Algebra i Logika \textbf{12}, 83--113 (1973) [Russian; English transl.\ Algebra and Logic \textbf{12}, 47--63
(1973)]

\bibitem{Ro80}
Rowen, L. H.: Polynomial Identities in Ring Theory. Academic Press, New York--London (1980)

\bibitem{sapirinherently}
Sapir, M. V.: Inherently nonfinitely based finite semigroups. Mat. Sb. \textbf{133}, 154--166 (1987) [Russian; English
transl. Math. USSR-Sb. \textbf{61}, 155--166 (1988)]

\bibitem{sapirburnside}
Sapir, M. V.: Problems of Burnside type and the finite basis property in varieties of semigroups. Izv.\ Akad.\ Nauk
SSSR, Ser.\ Mat. \textbf{51}, 319--340 (1987) [Russian; English transl.\ Math.\ USSR-Izv. \textbf{30}, 295--314 (1987)]

\bibitem{sapirinverse}
Sapir, M. V.: Identities of finite inverse semigroups. Internat.\ J. Algebra Comput. \textbf{3}, 115--124 (1993)

\bibitem{Sapir}
Sapir, M. V.: Combinatorics on Words with Applications. IBP-Litp 1995/32: Rapport de Recherche Litp, Universit\'e Paris
7 (1995). Available online under\\
\url{http://www.math.vanderbilt.edu/~msapir/ftp/course/course.pdf}

\bibitem{Serre}
Serre, J.-P.: Cours d'Arithmetique. Presses Universitaires de France, Paris (1980)

\bibitem{V}
Volkov, M. V.: On finite basedness of semigroup varieties. Mat. Zametki \textbf{45}, 12--23 (1989) [Russian; English
transl.\ Math.\ Notes \textbf{45}, 187--194 (1989)]

\bibitem{volkovjaponicae}
Volkov, M. V.: The finite basis problem for finite semigroups. Sci.\ Math.\ Jpn. \textbf{53}, 171--199 (2001)

\bibitem{Vo04}
Volkov, M. V.: Reflexive relations, extensive transformations and piecewise testable languages of a given height.
Internat. J. Algebra Comput. \textbf{14}, 817--827 (2004)

\bibitem{VoGo04}
Volkov, M. V., Goldberg, I. A.: The finite basis problem for monoids of triangular boolean matrices. In: Algebraic
Systems, Formal Languages, and Conventional and Unconventional Computation Theory, Surikaisekikenkyusho Kokyuroku
\textbf{1366}, pp. 205--214, Research Institute for Mathematical Sciences, Kyoto University, Kyoto (2004)
\end{thebibliography}
\end{document}